\documentclass[reqno,12pt]{amsart}
\usepackage{amsfonts}
\usepackage{amsmath}
\usepackage{amsthm}
\usepackage{url}

\newtheorem{theorem}{Theorem}[section]
\newtheorem{lemma}[theorem]{Lemma}
\newtheorem{proposition}[theorem]{Proposition}
\newtheorem{corollary}[theorem]{Corollary}
\newtheorem{definition}[theorem]{Definition}
\newtheorem{example}[theorem]{Example}
\newtheorem{remark}[theorem]{Remark}

\newtheorem*{theorem*}{Theorem}
\newtheorem*{lemma*}{Lemma}
\newtheorem*{proposition*}{Proposition}
\newtheorem*{corollary*}{Corollary}
\newtheorem*{definition*}{Definition}
\newtheorem*{example*}{Example}
\newtheorem*{remark*}{Remark}

\newcommand\uncorr[1]{\begin{corollary*}{\sl \  1}\end{corollary*}}

\numberwithin{equation}{section}

\allowdisplaybreaks[1]

\newcommand\ka{\kappa}

\title[Series and integrals associated with root systems]
{Hypergeometric and
basic hypergeometric series and integrals
associated with root systems}
\author{Michael J.\ Schlosser}
\address{Fakult\"at f\"ur Mathematik, Universit\"at Wien,
Oskar-Morgenstern-Platz~1, A-1090 Vienna, Austria}
\email{michael.schlosser@univie.ac.at}
\urladdr{http://www.mat.univie.ac.at/{\textasciitilde}schlosse}
\thanks{Partly supported by FWF Austrian Science Fund
grant F50-08.}

\begin{document}

\begin{abstract}
We give an overview of some of the main results from the theories of
hypergeometric and basic hypergeometric series and integrals
associated with root systems. In particular, we list a number of
summations, transformations and explicit evaluations for such
multiple series and integrals.
We concentrate on such results which do not directly extend to the
elliptic level.
This text is a provisional version of a chapter on
hypergeometric and basic hypergeometric series and integrals
associated with root systems for volume 2 of the new Askey--Bateman
project which deals with ``Multivariable special functions''. 
\end{abstract}

\maketitle
\tableofcontents

\section{Introduction}

Hypergeometric series associated with root systems
first appeared implicitly in the 1972 work of
Ali\v{s}auskas, Jucys and Jucys~\cite{AJJ} and
Chac\'on, Ciftan and Biedenharn~\cite{CCB}
in the context of the representation theory
of the unitary groups, more precisely,
as the multiplicity-free Wigner and Racah coefficients
($3j$ and $6j$-symbols) of the group $\mathrm{SU}(n+1)$.
A few years later, Holman, Biedenharn and Louck~\cite{HBL}
investigated these coefficients more explicitly as generalized
hypergeometric series and obtained a first summation theorem for these.
The series in question have explicit summands and contain
the Weyl denominator of the root system $\mathrm A_n$,
and can thus be considered as hypergeometric series associated with
this root system.
(These series are not to be confused with the hypergeometric functions
associated with root systems considered in Chapter~8 of this volume,
which generalize the spherical functions on noncompact Riemannian
symmetric spaces.)
Subsequently, $\mathrm A_n$
hypergeometric series were shown to satisfy various extensions
of well-known identities for classical hypergeometric series
\cite{G1,G2,H,M80}.
For example, using the usual Pochhammer symbol notation
for the shifted factorial (see \eqref{Poch}), 
Holman's~\cite{H} $\mathrm A_n$ extension of the
terminating balanced Pfaff--Saalsch\"utz ${}_3F_2$ summation is
\begin{align}
\sum_{k_1,\dots,k_n=0}^{N_1,\dots,N_n}&\Bigg(\prod_{1\le i<j\le n}
\frac{x_i+k_i-x_j-x_k}{x_i-x_j}\prod_{i,j=1}^n
\frac{(-N_j+x_i-x_j)_{k_i}}{(1+x_i-x_j)_{k_i}}\notag\\
&\qquad\times\prod_{i=1}^n\frac{(a+x_i)_{k_i}\,(b+x_i)_{k_i}}
{(c+x_i)_{k_i}\,(a+b-c+1-|N|+x_i)_{k_i}}\Bigg)\notag\\
&\qquad=\frac{(c-a)_{|N|}\,(c-b)_{|N|}}
{\prod_{i=1}^n(c+x_i)_{N_i}\,(c-a-b+|N|-N_i-x_i)_{N_i}},
\end{align} 
where, throughout this chapter,
$|N|:=N_1+\cdots+N_n$.

An extensive study of the \textit{basic} analogue, or $q$-analogue,
of $\mathrm A_n$ series was initiated by Milne in a series of papers
\cite{M1,M2,M3}.
The following application of the
``fundamental theorem of  $\mathrm A_n$ series'' \cite[Theorem~1.49]{M2},
\begin{equation}\label{fthmmilne}
\sum_{\substack{k_1,\dots,k_n\ge0\\|k|=N}}
\prod_{1\le i<j\le n}\frac{x_iq^{k_i}-x_jq^{k_j}}{x_i-x_j}
\prod_{i,j=1}^n\frac{(a_jx_i/x_j;q)_{k_i}}{(qx_i/x_j;q)_{k_i}}
=\frac{(a_1\cdots a_n;q)_N}{(q;q)_N},
\end{equation}
where we are using the usual $q$-Pochhammer symbol notation
for the $q$-shifted factorial (see \eqref{qPoch}),
demonstrates a phenomenon which is typical for the
$\mathrm A_n$ theory:
In \eqref{fthmmilne}, let $n\mapsto n+1$ and replace
$k_{n+1}$ by $N-(k_1+\dots+k_n)$. Then, after further replacing the
variables $a_i$ by $c_i$, for $i=1,\dots,n$, $a_{n+1}$ by $q^{-N}b/a$ and
$x_{n+1}$ by $q^{-N}/a$, respectively,
the following terminating $\mathrm A_n$ ${}_6\phi_5$ summation is obtained:
\begin{align}\label{an65eq}
\sum_{\substack{k_1,\dots,k_n\ge0\\|k|\le N}}\Bigg(&
\prod_{1\le i<j\le n}\frac{x_iq^{k_i}-x_jq^{k_j}}{x_i-x_j}
\prod_{i=1}^n\frac{1-ax_iq^{k_i+|k|}}{1-ax_i}
\prod_{i,j=1}^n\frac{(c_jx_i/x_j;q)_{k_i}}{(qx_i/x_j;q)_{k_i}}\notag\\&\times
\prod_{i=1}^n\frac{(ax_i;q)_{|k|}\,(bx_i;q)_{k_i}}
{(ax_iq/c_i;q)_{|k|}\,(ax_iq^{1+N};q)_{k_i}}\cdot
\frac{(q^{-N};q)_{|k|}}{(aq/b;q)_{|k|}}
\left(\frac{aq^{1+N}}{bc_1\cdots c_n}\right)^{|k|}
\Bigg)\notag\\
&\qquad\qquad\qquad\qquad\quad=\frac{(aq/bc_1\cdots c_n;q)_N}{(aq/b;q)_N}
\prod_{i=1}^n\frac{(ax_iq;q)_N}{(ax_iq/c_i;q)_N}.
\end{align}

By application of the one-variable $q$-binomial theorem,
it follows that another consequence of \eqref{fthmmilne}
is the following $\mathrm A_n$ extension of the nonterminating
$q$-binomial theorem:
\begin{equation}\label{ntqbin1}
\sum_{k_1,\dots,k_n\ge0}
\prod_{1\le i<j\le n}\frac{x_iq^{k_i}-x_jq^{k_j}}{x_i-x_j}
\prod_{i,j=1}^n\frac{(a_jx_i/x_j;q)_{k_i}}{(qx_i/x_j;q)_{k_i}}\cdot z^{|k|}
=\frac{(a_1\cdots a_nz;q)_\infty}{(z;q)_\infty},
\end{equation}
valid for $|q|<1$ and $|z|<1$.

While $\mathrm A_n$ (basic) hypergeometric series have also been
referred to as $\mathrm{SU}(n)$ or $\mathrm{U}(n)$ series,
the terminology \textit{(basic) hypergeometric series associated
to the root system $\mathrm A_n$}, or simply
\textit{$\mathrm A_n$ (basic) hypergeometric series},
is preferred by most authors nowadays.

A further important development was Gustafson's introduction
of very-well-poised series for other root systems \cite{G3,G4}.
Gustafson also introduced related multivariate integrals
associated with root systems \cite{DG,G5,G6,G7,G8}.
In this setting the multiple series or integrals are classified
according to the type of specific factors (such as a Weyl denominator)
appearing in the summand or integrand.

Most of the known results for multivariate (basic) hypergeometric 
series and integrals associated with root systems indeed concern
classical root systems, while only sporadically
summations or transformations for series or integrals associated with
exceptional root systems have been obtained
\cite{vdB,G5,I1,I2,IT}.

The root system classification appears to be very useful
and one would hope that the various relations satisfied by
the series or integrals can be interpreted in terms of root systems
or even Lie theory.
Although the type of series in question first arose
(in the limit $q\rightarrow 1$)
in the representation theory of compact Lie groups,
many questions remain open about this connection.
While a (quantum) group interpretation for the $\mathrm A_n$ type of series
has been given by Rosengren in \cite{R2} --- even for the elliptic extension
of the series, surveyed in Chapter~6 of this volume --- 
no analogous interpretations for the other root systems
have yet been revealed.

In many instances various types
(still referring to the root system classification)
of series/integrals can be combined with each other after
which one obtains series/integrals of some ``mixed type''
for which the correct classification is not really clear.
The conclusion is that the root system classification of the
series/integrals considered here is only rough and not always precise.

In terms of the rough classification in
\cite[Section~2]{B}, \cite[Section~1]{BS}, and \cite[Section~5]{M9},
a multivariate series
$\sum_{k_1,\dots,k_{n+1}}S_{k_1,\dots,k_{n+1}}$ is considered to be an
\textit{$\mathrm A_n$ hypergeometric series},
if the summand $S_{k_1,\dots,k_{n+1}}$ contains the factor
\begin{subequations}
\begin{equation}\label{wda}
\prod_{1\le i<j\le n+1}(x_i+k_i-x_j-k_j).
\end{equation}
It is considered to be an \textit{$\mathrm A_n$ basic
hypergeometric series}, if it contains the factor
\begin{equation}\label{wdaq}
\prod_{1\le i<j\le n+1}\big(x_iq^{k_i}-x_jq^{k_j}\big).
\end{equation}
\end{subequations}
If we take the sum over $k_1+\dots+k_{n+1}=N$, we may
replace $k_{n+1}$ by $N-|k|$ (where, as before, $|k|=k_1+\dots+k_n$)
and the two respective products in \eqref{wda}/\eqref{wdaq}
can be written as
\begin{subequations}
\begin{equation}\label{wda1}
\prod_{1\le i<j\le n}(x_i+k_i-x_j-k_j)\,\prod_{i=1}^n(a+x_i+k_i+|k|),
\end{equation}
where $x_{n+1}$ was substituted by $-a-N$, and
\begin{equation}\label{wdaq1}
\Big({-}aq^{|k|}\Big)^{-n}\prod_{1\le i<j\le n}\big(x_iq^{k_i}-x_jq^{k_j}\big)\,
\prod_{i=1}^n\big(1-ax_iq^{k_i+|k|}\big),
\end{equation}
\end{subequations}
where $x_{n+1}$ was substituted by $q^{-N}/a$,
respectively.

%Similar definitions apply to other root systems.
%In particular, a \textit{$\mathrm D_n$ hypergeometric series} contains
%the factor
%\begin{subequations}
%\begin{equation}
%\prod_{1\le j<k\le n}\Big((u_j+y_j-u_k-y_k)(u_j+y_j+u_k+y_k)\Big)
%\end{equation}
%and a \textit{$D_n$ basic hypergeometric series} the factor
%\begin{equation}
%\prod_{1\le j<k\le n}\Big((u_jq^{y_j}-u_kq^{y_k})(1-u_ju_kq^{y_j+y_k})\Big).
%\end{equation}
%\end{subequations}
Likewise, a \textit{$\mathrm C_n$ hypergeometric series} contains
the factor
\begin{subequations}
\begin{equation}
\prod_{1\le i<j\le n}(x_i+k_i-x_j-k_j)\prod_{1\le i\le j\le n}(x_i+k_i+x_j+k_j)
\end{equation}
and a \textit{$\mathrm C_n$ basic hypergeometric series} the factor
\begin{equation}
\prod_{1\le i<j\le n}\big(x_iq^{k_i}-x_jq^{k_j}\big)
\prod_{1\le i\le j\le n}\big(1-x_ix_jq^{k_i+k_j}\big).
\end{equation}
\end{subequations}
We omit giving similar definitions for other root systems.
The above factors may be associated with the Weyl denominators
$\prod_{\alpha>0}(1-e^{-\alpha})$ with the product taken over the
positive roots in the root system.
(Weyl denominators similarly appear in
Chapter~8, Section~4.2, of this volume.)

A very similar classification applies to hypergeometric
\textit{integrals} associated with root systems,
by considering certain factors of the integrand.
For specific examples, see Section~\ref{sec:int}.

%Sometimes the association with a specific root system is not clear.
%For instance, there are multivariate basic hypergeometric series
%which some authors \cite{B,BS,Schl0} have formally associated
%with the root system $\mathrm D_n$ but which other authors \cite{SW}
%prefer to associate with $\mathrm A_n$.
%Since the $\mathrm D_n$ series considered in \cite{B,BS,Schl0} actually
%arise by combining $\mathrm A_n$ with $\mathrm C_n$ series,
%and they share features of $\mathrm A_n$ and $\mathrm C_n$ series,
%one might just refer to them as of ``mixed root system type''
%(or short, of ``mixed type'').

A special feature of the theory of hypergeometric or basic hypergeometric
functions associated with root systems
is that often there exist several different identities
for one and the same root system that extend a particular
one-variable identity. See the various
$\mathrm A_n$ ${}_3\phi_2$ and $\mathrm A_n$ ${}_2\phi_1$ summations
given by Milne in \cite{M6}. 
At this point we would also like to mention that various special
$\mathrm A_n$ hypergeometric series
possess rich structures of symmetry; these were made explicit
by Kajihara~\cite{K3}.

This chapter deals to a large extent with the multivariate \textit{basic}
hypergeometric theory.
The reason is that most results for ordinary hypergeometric series
have basic hypergeometric analogues\footnote{Identities for
basic hypergeometric series
reduce to those for ordinary hypergeometric series by taking the limit
$q\to 1$ in an appropriate manner~\cite[Section 1.2]{GR}, see for
instance the derivation of \eqref{an2h2cglN} from \eqref{an1psi1cglN}.},
so it does not make sense to treat the hypergeometric case
separately.

A good amount of the theory of basic hypergeometric series
associated with root systems has recently been generalized to the
\textit{elliptic} level.
%For an introduction to one-variable
%elliptic hypergeometric functions, see Volume~3.
For an introduction to elliptic hypergeometric functions
associated with root systems, see Chapter~6 of this volume.
In the present chapter we emphasize some
general facts about series associated with root systems which
are important for understanding the nature of the series,
but otherwise, to avoid overlap, mainly
focus on parts of the theory which do \textit{not}
directly extend to the elliptic level.
%This in particular concerns some
%identities involving only very few parameters (which can be obtained
%from more general identities by taking specific limits,
%say by letting some parameters
%tend to zero or infinity, a simple technique which is not available
%when dealing with elliptic hypergeometric series), and also
%infinite (i.e.\ nonterminating, possibly even multilateral) series.
This in particular concerns
identities obtained as confluent limits of more general identities
and identities for nonterminating and/or multilateral series.

The following sections are devoted to multivariate identities,
ranging from very simple identities to more complicated ones.
In particular, in Sections~\ref{sec:2} and \ref{sec:int}
various summations, transformations, and integral evaluations
are reviewed. Section~\ref{sec:macd} surveys the theory of
basic hypergeometric series with Macdonald polynomial argument.
The chapter concludes with Section~\ref{sec:appl},
containing a brief discussion on applications of
basic hypergeometric series associated with root systems.
\smallskip

\noindent
\textbf{Acknowledgements:} I would like to thank Gaurav Bhatnagar,
Tom Koornwinder, Stephen Milne, Hjalmar Rosengren, Jasper Stokman and
Ole Warnaar for careful reading and many valuable comments.
The author's research was partially supported by Austrian Science Fund
grant F50-08.

\section{Some identities for (basic) hypergeometric series
associated with root systems}\label{sec:2}

A large number of identities for hypergeometric and basic hypergeometric
series associated with root systems has appeared in the literature.
Due to space limitations, we only provide a small representative
selection of identities. Nevertheless, they are meant to give a flavor
of the expressions which typically occur in the multivariate theory.
For more details the reader is pointed to specific literature.

We use the following notations for shifted and $q$-shifted factorials
(which are also referred to as Pochhammer and
$q$-Pochhammer symbols, respectively):
\begin{equation}\label{Poch}
(a)_k=\begin{cases}1&k=0\\
a(a+1)\dots(a+k-1)&k=1,2,\dots,\\
\big((a-1)(a-2)\dots(a+k)\big)^{-1}&k=-1,-2,\dots,
\end{cases}
\end{equation}
\begin{subequations}\label{qPoch}
\begin{align}
(a;q)_k&=\begin{cases}1&k=0\\
(1-a)(1-aq)\dots(1-aq^{k-1})&k=1,2,\dots,\\
\big((1-aq^{-1})(1-aq^{-2})\dots(1-aq^k)\big)^{-1}&k=-1,-2,\dots,
\end{cases}\\
\intertext{and}
(a;q)_\infty&=\prod_{i\ge 0}(1-aq^i).
\end{align}
\end{subequations}
When dealing with products of shifted and $q$-shifted factorials,
we frequently use the shorthand notations
\begin{equation*}
(a_1,\dots,a_n)_j=(a_1)_j\cdots(a_n)_j\quad\text{and}\quad
(a_1,\dots,a_n;q)_k=(a_1;q)_k\cdots(a_n;q)_k,
\end{equation*}
where $j$ is an integer, and $k$ is an integer or $\infty$.
 
This chapter reviews results for multivariate extensions associated with
root systems of the following univariate series,
whose definitions we give for self-containment.

Hypergeometric ${}_rF_s$ series and bilateral hypergeometric
${}_rH_s$ series are defined as
\begin{subequations}
\begin{align}
{}_rF_s\!\left[\begin{matrix}a_1,\dots,a_r\\
b_1,\dots,b_s\end{matrix};z\right]&=
\sum_{k=0}^\infty\frac{(a_1,\dots,a_r)_k}{(1,b_1,\dots,b_s)_k}z^k,\\
{}_rH_s\!\left[\begin{matrix}a_1,\dots,a_r\\
b_1,\dots,b_s\end{matrix};z\right]&=
\sum_{k=-\infty}^\infty\frac{(a_1,\dots,a_r)_k}{(b_1,\dots,b_s)_k}z^k.
\end{align}
\end{subequations}
Similarly,
basic hypergeometric ${}_r\phi_s$ series and bilateral basic
hypergeometric ${}_r\psi_s$ series are defined as
\begin{subequations}
\begin{align}
{}_r\phi_s\!\left[\begin{matrix}a_1,\dots,a_r\\
b_1,\dots,b_s\end{matrix};q,z\right]&=
\sum_{k=0}^\infty\frac{(a_1,\dots,a_r;q)_k}{(q,b_1,\dots,b_s;q)_k}
\left((-1)^kq^{\binom k2}\right)^{1+s-r}z^k,\\
{}_r\psi_s\!\left[\begin{matrix}a_1,\dots,a_r\\
b_1,\dots,b_s\end{matrix};q,z\right]&=
\sum_{k=-\infty}^\infty\frac{(a_1,\dots,a_r;q)_k}{(b_1,\dots,b_s;q)_k}
\left((-1)^kq^{\binom k2}\right)^{s-r}z^k.
\end{align}
\end{subequations}
See Slater's book \cite{Sl0} and
Gasper and Rahman's book \cite{GR} for the conditions of theses series
to terminate, to be balanced, and to be (very-)well-poised,
and for various identities satisfied by these series.

\subsection{Some useful elementary facts}\label{subsecel}

We start with a few elementary ingredients
which are useful for manipulating basic hypergeometric series
associated with root systems.

\begin{itemize}
\item[(i)]
A fundamental ingredient (for inductive proofs and
functional equations, etc.) is the following
partial fraction decomposition \cite[Section~7]{M5}:
\begin{equation}\label{pfd}
\prod_{i=1}^n\frac {1-tx_iy_i} {1-tx_i}=y_1y_2\cdots y_n+
\sum_{k=1}^n\frac {\prod_{i=1}^n(1-y_ix_i/x_k)}
{(1-tx_k)\,\prod_{\substack{i=1\\i\neq k}}^n(1-x_i/x_k)}.
\end{equation}
In particular, this identity can be used to prove the fundamental
theorem of $\mathrm A_n$ series in \eqref{fthmmilne}. 

A slightly more general partial fraction decomposition was
derived in \cite[Lemma~3.2]{Schl6}.
The identity in \eqref{pfd} can be obtained as a limiting case
of an elliptic
partial fraction decomposition of type $\mathrm A$,
cf.\ \cite[p.\ 451, Example~3]{WW}.
A related partial fraction decomposition of type $\mathrm D$ was established
in \cite[Lemma 4.14]{G2}.

\item[(ii)]
For simplifying products the following identity is useful
\cite[Lemma~4.3]{M6}:
\begin{equation}
\prod_{1\le i<j\le n}\frac{x_iq^{k_i}-x_jq^{k_j}}{x_iq^{m_i}-x_jq^{m_j}}
\prod_{i,j=1}^n\frac{(q^{m_i-k_j}x_i/x_j;q)_{k_i-m_i}}{(q^{1+m_i-m_j}x_i/x_j;q)_{k_i-m_i}}
=(-1)^{|k|-|m|}q^{-\binom{|k|-|m|+1}2}.
\end{equation}

\item[(iii)]
When reversing the order of the summations
\begin{equation*}
\sum_{k_1,\dots,k_n=0}^{N_1,\dots,N_n}S_{k_1,\dots,k_n}=
\sum_{k_1,\dots,k_n=0}^{N_1,\dots,N_n}S_{N_1-k_1,\dots,N_n-k_n},
\end{equation*}
it is convenient to use the fact that the following variant of an
``$\mathrm A_n$ $q$-binomial coefficient''
\begin{equation}\label{anqbinc}
\prod_{i,j=1}^n\frac{(qx_i/x_j;q)_{N_i}}
{(qx_i/x_j;q)_{k_i}\,(q^{1+k_i-k_j}x_i/x_j;q)_{N_i-k_i}}
\end{equation}
(the usual $q$-binomial coefficient is given in \eqref{qbinc})
remains unchanged after performing the simultaneous substitutions
$k_i\mapsto N_i-k_i$ and $x_i\mapsto q^{-N_i}/x_i$, for $i=1,\dots,n$
\cite[Remark~B.3]{Schl2}. 
\end{itemize}

\subsection{Some terminating
%and nonterminating
$\mathrm A_n$ $q$-binomial theorems}

The terminating $q$-binomial theorem can be written in the form
(cf.\ \cite[Ex.~1.2(vi)]{GR})
\begin{equation}\label{tqbin}
(z;q)_N=\sum_{k=0}^N\begin{bmatrix}N\\k\end{bmatrix}_q
(-1)^kq^{\binom k2}z^k,
\end{equation}
where
\begin{equation}\label{qbinc}
\begin{bmatrix}N\\k\end{bmatrix}_q=
\frac{(q;q)_N}{(q;q)_k(q;q)_{N-k}}
\end{equation}
is the $q$-binomial coefficient. (Here, $N$ denotes a nonnegative integer.)
In basic hypergeometric notation, this identity corresponds to a terminating
${}_1\phi_0$ summation. It can be immediately obtained from the
nonterminating $q$-binomial theorem (or ${}_1\phi_0$ summation,
cf.\ \cite[Equation (II.3)]{GR}),
\begin{equation}\label{ntqbin0}
{}_1\phi_0\!\left[\begin{matrix}a\\
-\end{matrix};q,z\right]=\frac{(az;q)_\infty}{(z;q)_\infty},
\end{equation}
valid for $|q|<1$ and $|z|<1$. To obtain \eqref{tqbin} from \eqref{ntqbin0},
replace $a$ and $z$ by $q^{-n}$ and $zq^n$, respectively.

We have the following three multisum identities (for equivalent forms,
where the Vandermonde determinant of type $\mathrm A$ \eqref{wdaq}
explicitly appears in the summands, 
see \cite[Theorems~5.44, 5.46, and 5.48]{M6}),
each involving the $\mathrm A_n$ $q$-binomial coefficient in
\eqref{anqbinc}:
\begin{subequations}
\begin{equation}\label{qbin1}
(z;q)_{|N|}=\sum_{k_1,\dots,k_n=0}^{N_1,\dots,N_n}
\prod_{i,j=1}^n\frac{(qx_i/x_j;q)_{N_i}}
{(qx_i/x_j;q)_{k_i}\,(q^{1+k_i-k_j}x_i/x_j;q)_{N_i-k_i}}\cdot
(-1)^{|k|}q^{\binom{|k|}2}z^{|k|},
\end{equation}
\begin{align}\label{qbin2}
\prod_{i=1}^n(zx_i;q)_{N_i}=\sum_{k_1,\dots,k_n=0}^{N_1,\dots,N_n}&\Bigg(
\prod_{i,j=1}^n\frac{(qx_i/x_j;q)_{N_i}}
{(qx_i/x_j;q)_{k_i}\,(q^{1+k_i-k_j}x_i/x_j;q)_{N_i-k_i}}\notag\\&\ \times
(-1)^{|k|}q^{\sum_{i=1}^n\binom{k_i}2}z^{|k|}\prod_{i=1}^nx_i^{k_i}\Bigg),
\end{align}
\begin{align}\label{qbin3}
\prod_{i=1}^n(zq^{|N|-N_i}/x_i;q)_{N_i}=\sum_{k_1,\dots,k_n=0}^{N_1,\dots,N_n}&\Bigg(
\prod_{i,j=1}^n\frac{(qx_i/x_j;q)_{N_i}}
{(qx_i/x_j;q)_{k_i}\,(q^{1+k_i-k_j}x_i/x_j;q)_{N_i-k_i}}\notag\\&\ \times
(-1)^{|k|}q^{\binom{|k|}2+\sum_{1\le i<j\le n}k_ik_j}z^{|k|}\prod_{i=1}^nx_i^{-k_i}\Bigg).
\end{align}
%where $e_2(k)$ is the second elementary symmetric function of
%$k_1,\dots,k_n$.

The summations in \eqref{qbin2} and \eqref{qbin3} are related by
reversal of sums as explained in Subsection~\ref{subsecel}(iii).
In the last two identities the variable $z$ is redundant.
Inclusion turns both sides into polynomials in $z$.

Here are two other terminating $\mathrm A_n$ $q$-binomial theorems
\cite[Theorems~5.52 and 5.50]{M6}:
\begin{equation}\label{qbinm1}
(z;q)_{N}=\sum_{\substack{k_1,\dots,k_n\ge0\\|k|\le N}}
\prod_{1\le i<j\le n}\frac{x_iq^{k_i}-x_jq^{k_j}}{x_i-x_j}
\prod_{i,j=1}^n(qx_i/x_j;q)_{k_i}^{-1}\cdot(q^{-N};q)_{|k|}\,q^{N|k|}z^{|k|},
\end{equation}
\begin{align}\label{qbinm2}
(z;q)_{N}=\sum_{\substack{k_1,\dots,k_n\ge0\\|k|\le N}}&\Bigg(
\prod_{1\le i<j\le n}\frac{x_iq^{k_i}-x_jq^{k_j}}{x_i-x_j}
\prod_{i,j=1}^n(qx_i/x_j;q)_{k_i}^{-1}\cdot(q^{-N};q)_{|k|}\notag\\&\ \times
(-1)^{(n-1)|k|}
q^{N|k|-\binom{|k|}2+n\sum_{i=1}^n\binom{k_i}2}
z^{|k|}\prod_{i=1}^nx_i^{nk_i-|k|}\Bigg).
\end{align}
Note that Equations \eqref{qbinm1} and \eqref{qbinm2} are equivalent
with respect to inverting the base $q\to q^{-1}$.

Yet another terminating $\mathrm A_n$ $q$-binomial theorem, implicit from
in \cite{BS}, is
\begin{align}\label{qbin4}
\prod_{i=1}^n(z/x_i;q)_{N}=\sum_{\substack{k_1,\dots,k_n\ge 0\\|k|\le N}}
&\Bigg(\prod_{1\le i<j\le n}\frac{x_iq^{k_i}-x_jq^{k_j}}{x_i-x_j}
\prod_{i,j=1}^n(qx_i/x_j;q)_{k_i}^{-1}\cdot(q^{-N};q)_{|k|}
\notag\\&\ \times
q^{N|k|+\sum_{1\le i<j\le n}k_ik_j}z^{|k|}
\prod_{i=1}^nx_i^{-k_i}(z/x_i;q)_{|k|-k_i}\Bigg).
\end{align}
\end{subequations}

Four of the terminating $\mathrm A_n$ $q$-binomial theorems of this subsection
are special cases of nonterminating $\mathrm A_n$ $q$-binomial theorems.
In particular, \eqref{qbin1} is a special case of \eqref{ntqbin1}.
Similarly, \eqref{qbin2} is a special case of \cite[Theorem~4.7]{LM}.
Further, \eqref{qbin4} is a special case of \cite[Theorem~5.42]{M6}
(which is the $b_1=\cdots=b_n=q$ special case of \eqref{1psi1b}).
Finally, \eqref{qbinm2} is a special case of \cite[Theorem~5.19]{BS}.

\subsection{Other terminating summations}

All the $\mathrm A_n$ $q$-binomial theorems (i.e., terminating
${}_1\phi_0$ summations)
listed in the previous subsection (and others not listed here,
as those implicit in \cite[Section~7.7]{RS}) admit generalizations to
summations involving more parameters.
These include in particular various
$\mathrm A_n$ ${}_2\phi_1$, ${}_3\phi_2$, ${}_6\phi_5$, or ${}_8\phi_7$
summations (see e.g.\ \cite{M6,RS,Schl5}).

The terminating balanced ${}_3\phi_2$ summation
(or $q$-Pfaff--Saalsch\"utz summation) is (cf.\ \cite[Equation (II.12)]{GR})
\begin{equation}
{}_3\phi_2\!\left[\begin{matrix}a,b,q^{-N}\\
c,abq^{1-N}/c\end{matrix};q,q\right]=\frac{(c/a,c/b;q)_N}{(c,c/ab;q)_N}.
\end{equation}
Here are two $\mathrm A_n$ ${}_3\phi_2$ summations (cf.\ \cite{M6}):
\begin{multline}\label{an32feq}
\sum_{k_1,\dots,k_n=0}^{N_1,\dots,N_n}
\prod_{1\le i<j\le n}\frac{x_iq^{k_i}-x_jq^{k_j}}{x_i-x_j}
\prod_{i,j=1}^n\frac{(q^{-N_j}x_i/x_j;q)_{k_i}}{(qx_i/x_j;q)_{k_i}}
\prod_{i=1}^n\frac{(ax_i;q)_{k_i}}
{(cx_i;q)_{k_i}}\cdot\frac{(b;q)_{|k|}}{(abq^{1-|N|}/c;q)_{|k|}}q^{|k|}\\
=\frac{(c/a;q)_{|N|}}{(c/ab;q)_{|N|}}
\prod_{i=1}^n\frac{(cx_i/b;q)_{N_i}}{(cx_i;q)_{N_i}},
\end{multline}
\begin{multline}\label{an32teq}
\sum_{k_1,\dots,k_n=0}^{N_1,\dots,N_n}
\prod_{1\le i<j\le n}\frac{x_iq^{k_i}-x_jq^{k_j}}{x_i-x_j}
\prod_{i,j=1}^n\frac{(q^{-N_j}x_i/x_j;q)_{k_i}}{(qx_i/x_j;q)_{k_i}}
\prod_{i=1}^n\frac{(ax_i,bx_i;q)_{k_i}}
{(cx_i,abx_iq^{1-|N|}/c;q)_{k_i}}q^{|k|}\\
=(c/a,c/b;q)_{|N|}
\prod_{i=1}^n(cx_i,cx_iq^{|N|-N_i}/ab;q)_{N_i}^{-1}.
\end{multline}

Several other $\mathrm A_n$ ${}_3\phi_2$ summations are given in \cite{M6}.
For instance, a simple polynomial argument applied to \eqref{an32feq} yields
\begin{multline}\label{an32fdeq}
\sum_{\substack{k_1,\dots,k_n\ge0\\|k|\le N}}
\prod_{1\le i<j\le n}\frac{x_iq^{k_i}-x_jq^{k_j}}{x_i-x_j}
\prod_{i,j=1}^n\frac{(a_jx_i/x_j;q)_{k_i}}{(qx_i/x_j;q)_{k_i}}
\prod_{i=1}^n\frac{(bx_i;q)_{k_i}}{(cx_i;q)_{k_i}}\cdot
\frac{(q^{-N};q)_{|k|}}{(a_1\cdots a_nbq^{1-N}/c;q)_{|k|}}q^{|k|}\\
=\frac{(c/b;q)_N}{(c/a_1\cdots a_nb;q)_N}
\prod_{i=1}^n\frac{(cx_i/a_i;q)_N}{(cx_i;q)_N}.
\end{multline}

Here is another terminating balanced ${}_3\phi_2$ summation
(\cite[Theorem~4.3]{ML} rewritten)
which may be considered of ``mixed-type'':
\begin{align}\label{dn32fdeq}
\sum_{\substack{k_1,\dots,k_n\ge0\\|k|\le N}}
\Bigg(&\prod_{1\le i<j\le n}\frac{x_iq^{k_i}-x_jq^{k_j}}{x_i-x_j}
\frac 1{(x_ix_j;q)_{k_i+k_j}}
\prod_{i,j=1}^n\frac{(a_jx_i/x_j,x_ix_j/a_j;q)_{k_i}}
{(qx_i/x_j;q)_{k_i}}\notag\\&\times
\frac{(q^{-N};q)_{|k|}}{\prod_{i=1}^n(bx_iq^{-N},qx_i/b;q)_{k_i}}q^{|k|}\Bigg)
=\prod_{i=1}^n\frac{(qa_i/bx_i,qx_i/a_ib;q)_N}{(q/bx_i,qx_i/b;q)_N}.
\end{align}
Again, a simple polynomial argument can be applied to transform this
summation to another one, in this case to a sum over a rectangular region 
(see \cite[Theorem~4.2]{ML}), which we do not state here explicitly. 

Among the most general summations for basic hypergeometric series
associated with root systems are various multivariate
${}_8\phi_7$ Jackson summations.
In the univariate case, Jackson's terminating balanced very-well-poised
${}_8\phi_7$ summation is (cf.\ \cite[Equation (II.22)]{GR})
\begin{equation}
{}_8\phi_7\!\left[\begin{matrix}a,\,qa^{\frac 12},-qa^{\frac 12},b,c,d,e,q^{-N}\\
a^{\frac 12},-a^{\frac 12},aq/b,aq/c,aq/d,aq/e,aq^{N+1}\end{matrix};q,q\right]
=\frac{(aq,aq/bc,aq/bd,aq/cd;q)_N}{(aq/b,aq/c,aq/d,aq/bcd;q)_N},
\end{equation}
where $a^2q=bcdeq^{-N}$.
Some of the multivariate Jackson summations
have been extended to the level
of elliptic hypergeometric series and are partly covered in
Section~6.3 of this volume.
One of the most important is the following $\mathrm A_n$ Jackson summation
\cite[Theorem~6.14]{M7}:
\begin{multline}\label{an87eq}
\sum_{k_1,\dots,k_n=0}^{N_1,\dots,N_n}\Bigg(
\prod_{1\le i<j\le n}\frac{x_iq^{k_i}-x_jq^{k_j}}{x_i-x_j}
\prod_{i=1}^n\frac{1-ax_iq^{k_i+|k|}}{1-ax_i}
\prod_{i,j=1}^n\frac{(q^{-N_j}x_i/x_j;q)_{k_i}}{(qx_i/x_j;q)_{k_i}}\\\times
\prod_{i=1}^n\frac{(ax_i;q)_{|k|}\,(dx_i,a^2x_iq^{1+N_i}/bcd;q)_{k_i}}
{(ax_iq^{1+|N|};q)_{|k|}\,(ax_iq/b,ax_iq/c;q)_{k_i}}\cdot
\frac{(b,c;q)_{|k|}}{(aq/d,bcdq^{-|N|}/a;q)_{|k|}}
q^{|k|}\Bigg)\\
=\frac{(aq/bd,aq/cd;q)_{|N|}}{(aq/d,aq/bcd;q)_{|N|}}
\prod_{i=1}^n\frac{(ax_iq,ax_iq/bc;q)_{N_i}}
{(ax_iq/b,ax_iq/c;q)_{N_i}}.
\end{multline}
It was initially proved by partial fraction decompositions and functional
equations; a more direct proof (which extends to the elliptic level)
utilizes partial fraction decompositions and induction \cite{R1}.

From \eqref{an87eq}, by multivariable matrix inversion,
the following $\mathrm A_n$ Jackson summation was deduced in
\cite[Theorem~4.1]{Schl5}:
\begin{multline}\label{an87eq2}
\sum_{k_1,\dots,k_n=0}^{N_1,\dots,N_n}\Bigg(
\prod_{1\le i<j\le n}\frac{x_iq^{k_i}-x_jq^{k_j}}{x_i-x_j}
\prod_{i=1}^n\frac{(bcd/ax_i;q)_{|k|-k_i}\,(d/x_i;q)_{|k|}\,
(a^2x_iq^{1+|N|}/bcd;q)_{k_i}}
{(d/x_i;q)_{|k|-k_i}\,(bcdq^{-N_i}/ax_i;q)_{|k|}\,
(ax_iq/d;q)_{k_i}}\\\times
\prod_{i,j=1}^n\frac{(q^{-N_j}x_i/x_j;q)_{k_i}}{(qx_i/x_j;q)_{k_i}}\cdot
\frac{(1-aq^{2|k|})}{(1-a)}
\frac{(a,b,c;q)_{|k|}}
{(aq^{1+|N|},aq/b,aq/c;q)_{|k|}}
q^{|k|}\Bigg)\\
=\frac{(aq,aq/bc;q)_{|N|}}{(aq/b,aq/c;q)_{|N|}}
\prod_{i=1}^n\frac{(ax_iq/bd,ax_iq/bc;q)_{N_i}}
{(ax_iq/d,ax_iq/bcd;q)_{N_i}}.
\end{multline}
(Its elliptic extension is deduced in \cite{RS2}.)
Both summations, \eqref{an87eq} and \eqref{an87eq2}, which are summed
over rectangular regions, can be turned to summations over a
tetrahedral region (or simplex) $\{k_1,\dots,k_n\ge0,\ |k|\le N\}$
by a polynomial argument, a standard procedure in the multivariate theory,
used extensively in \cite{M6}.
Other Jackson summations which have been extended to the elliptic level
are the $\mathrm C_n$ Jackson summations in \cite[Theorem~4.1]{DG}
(independently derived in \cite[Theorem~6.13]{ML})
and \cite[Theorem~4.3]{Schl1},
the $\mathrm D_n$ Jackson summations (also referred to as
$\mathrm A_n$ Jackson summations by some authors)
of \cite[Theorem~7]{B} and \cite[Theorem~5.6]{Schl0},
the $\mathrm A_n$ and $\mathrm D_n$ Jackson summations in
\cite[Section~11]{BS2},
and the $\mathrm{BC}_n$ Jackson summation in \cite[Theorem~3]{vDS}
(also derived in \cite[$p\to 0$ in Theorem~2.1]{R0})
which was originally conjectured in \cite{W1}.

The following $\mathrm A_n$ Jackson summation is due to
Gustafson and Rakha \cite[Theorem~1.2]{GuR}
(but stated here as in \cite{R5}
where it has been extended to the elliptic level):
\begin{multline}\label{grsum}
\sum_{\substack{k_1,\dots,k_n\ge 0\\|k|\le N}}\Bigg(
\prod_{1\le i<j\le n}\frac{x_iq^{k_i}-x_jq^{k_j}}{x_i-x_j}\,
(x_ix_j;q)_{k_i+k_j}
\prod_{1\le i,j\le n}(qx_i/x_j;q)_{k_i}^{-1}\,
\prod_{i=1}^n\frac{1-ax_iq^{k_i+|k|}}{1-ax_i}\\\times
\prod_{i=1}^n\frac{(ax_i;q)_{|k|}}{(aq/x_i;q)_{|k|-k_i}\,(ax_iq^{1+N_i};q)_{k_i}}
\cdot\prod_{j=1}^4\frac{\prod_{i=1}^n(x_ib_j;q)_{k_i}}{(aq/b_j;q)_{|k|}}
\cdot(q^{-N};q)_{|k|}\, q^{|k|}\Bigg)\\
=\big(aq/b_1,aq/b_2,aq/b_3,aq/b_1b_2b_3X^2;q\big)_N^{-1}\;
\prod_{i=1}^n\frac{(ax_iq;q)_N}{(aq/x_i;q)_N}\\\times
\begin{cases}
\big(aq/X,aq/b_1b_2X,aq/b_1b_3X,aq/b_2b_3X;q\big)_N,&
\quad \text{if $n$ is odd},\\
\big(aq/b_1X,aq/b_2X,aq/b_3X,aq/b_1b_2b_3X;q\big)_N,&
\quad \text{if $n$ is even},
\end{cases}
\end{multline}
where $X=x_1\cdots x_n$,
under the assumption that $a^2q^{N+1}=b_1b_2b_3b_4X^2$.

Two similar multivariate ${}_8\phi_7$ summations
(at the elliptic level) are established in \cite{R5}, and two others
are conjectured in \cite[Conjectures~6.2 and 6.5]{SW}.
The latter actually look more complicated, the sums running over pairs
of partitions whose Ferrers diagrams differ by a horizontal strip
(cf.\ \cite{Mac0}). Those summations may play a role in the construction
of Askey--Wilson polynomials of type $\mathrm A$.

An $\mathrm A_n$ Jackson summation of a quite different type, intimately
related to Macdonald polynomials \cite{Mac0},
has been derived in \cite[Theorem~4.1]{Schl4}:
\begin{multline}\label{an87neq}
\sum_{\substack{k_1,\dots,k_n\ge 0\\|k|\le N}}\Bigg(
\prod_{i,j=1}^n\frac{{(qx_i/t_ix_j;q)}_{k_i}}
{{(qx_i/x_j;q)}_{k_i}}\,
\prod_{1\le i<j\le n}\frac{(t_jx_i/x_j;q)_{k_i-k_j}}
{(qx_i/t_ix_j;q)_{k_i-k_j}}\,\frac 1{x_i-x_j}\\\times
\det_{1\le i,j \le n}
\left[(x_iq^{k_i})^{n-j}
\left(1-t_i^{j-n-1}\frac{1-t_0x_iq^{k_i}}
{1-t_0x_iq^{k_i}/t_i}
\prod_{s=1}^n \frac{x_iq^{k_i}-x_s}
{x_iq^{k_i}/t_i-x_s}\right)\right]\\[.3em]\times
\prod_{i=1}^n
\frac{(dq^{-N}/t_0x_i;q)_{|k|-k_i}\,
(t_0x_iq/t_i,bx_i,t_0^2x_iq^{1+N}/bdt_1\cdots t_n;q)_{k_i}}
{(dt_iq^{-N}/t_0x_i;q)_{|k|-k_i}\,
(t_0x_iq,t_0x_iq/dt_i,t_0x_iq^{1+N}/t_i;q)_{k_i}}\\\times
\frac{(d,q^{-N};q)_{|k|}}
{(bdq^{-N}/t_0,t_0q/bt_1\cdots t_n;q)_{|k|}}\,
q^{\sum_{i=1}^n(2-i)k_i}\prod_{i=1}^nt_i^{(i-1)k_i+\sum_{j=i+1}^nk_j}\Bigg)\\=
\frac{(t_0q/b,t_0q/bdt_1\cdots t_n;q)_N}
{(t_0q/bd,t_0q/bt_1\cdots t_n;q)_N}
\prod_{i=1}^n\frac{(t_0x_iq/t_i,t_0x_iq/d;q)_N}
{(t_0x_iq,t_0x_iq/dt_i;q)_N}.
\end{multline}

A similar Jackson sum of type $\mathrm C_n$ has been conjectured in
\cite[Conjecture~4.5]{Schl4}. For the $\mathrm A_n$ identity in
\eqref{an87neq} and the similar $\mathrm C_n$ identity from \cite{Schl4}
(conjectured) elliptic extensions have not yet been established. The
difficulty stems from the special type of determinants
(which do not allow termwise elliptic extension) appearing in the
respective summands of the series.

The $\mathrm A_n$ Jackson summation in \eqref{an87neq} is also
remarkable in the sense that no corresponding multivariate Bailey
transformation has yet been found or conjectured. (The same applies to
\cite[Conjecture~4.5]{Schl4}.)
The other Jackson summations which we have discussed in this subsection
all can be generalized to transformations. Since in this chapter we are
mainly concerned with identities that do not directly extend to the
elliptic setting
%(while we are still listing some Jacksons summations due
%to their fundamental importance)
we are not reproducing any of the
multivariate Bailey transformations here.
(The only exception is the $\mathrm C_n$ nonterminating Bailey
transformation in \eqref{cnnt109gl}, as nonterminating series do
not admit a direct elliptic extension, for the reason of convergence.)
For a discussion on multivariate extensions of
Bailey's ${}_{10}\phi_9$ transformation, see Chapter~6,
Subsections~6.3.3 and 6.3.4 of this volume.

\subsection{Some multilateral summations}

Dougall's bilateral ${}_2H_2$ summation~\cite[Section~13]{D} is
\begin{equation}
{}_2H_2\!\left[\begin{matrix}a,b\\c,d\end{matrix};1\right]=
\frac{\Gamma(1-a)\Gamma(1-b)\Gamma(c)\Gamma(d)\Gamma(c+d-a-b-1)}
{\Gamma(c-a)\Gamma(c-b)\Gamma(d-a)\Gamma(d-b)},
\end{equation}
where the series either terminates, or $\mathrm{Re}(c+d-a-b-1)>0$
for convergence.
This identity does not admit a direct basic extension
(as a closed form ${}_2\psi_2$ summation with general parameters does not
exist). A related, similar looking identity is
Ramanujan's ${}_1\psi_1$ summation theorem (cf.\ \cite[Equation~(II.29)]{GR}),
\begin{equation}
{}_1\psi_1\!\left[\begin{matrix}a\\b\end{matrix};q,z\right]=
\frac{(q,b/a,az,q/az;q)_\infty}{(b,q/a,z,b/az;q)_\infty},
\end{equation}
where $|b/a|<|z|<1$.

An $\mathrm A_n$ extension of Dougall's ${}_2H_2$ summation theorem was
proved by Gustafson in \cite[Theorem~1.11]{G2} by induction and residue
calculus:
\begin{multline}
\sum_{k_1,\dots,k_n=-\infty}^\infty
\prod_{1\le i<j\le n}\frac{x_i+k_i-x_j-k_j}{x_i-x_j}
\prod_{i=1}^n\prod_{j=1}^{n+1}\frac{(a_j+x_i)_{k_i}}{(b_j+x_i)_{k_i}}\\
=\frac{\Gamma\Big({-}n+\sum_{j=1}^{n+1}(b_j-a_j)\Big)
\prod_{i=1}^n\prod_{j=1}^{n+1}\Gamma(1-a_j-x_i)\,\Gamma(b_j+u_i)}
{\prod_{i,j=1}^{n+1}\Gamma(b_j-a_i)
\prod_{1\le i<j\le n}\Gamma(1-x_i+x_j)\,\Gamma(1+x_i-x_j)},
\end{multline}
provided $\mathrm{Re}\Big(\sum_{j=1}^{n+1}(b_j-a_j)\Big)>n$.

Similarly, an $\mathrm A_n$ extension of Ramanujan's ${}_1\psi_1$
summation theorem was proved in \cite[Theorem~1.17]{G2}:
\begin{multline}\label{1psi1gus}
\sum_{k_1,\dots,k_n=-\infty}^{\infty}
\prod_{1\le i<j\le n}\frac{x_iq^{k_i}-x_jq^{k_j}}{x_i-x_j}
\prod_{i,j=1}^n\frac{(a_jx_i/x_j;q)_{k_i}}{(b_jx_i/x_j;q)_{k_i}}\cdot
z^{|k|}\\
=\frac{(a_1\cdots a_nz,q/a_1\cdots a_nz;q)_{\infty}}
{(z,b_1\cdots b_nq^{1-n}/a_1\cdots a_nz;q)_{\infty}}
\prod_{i,j=1}^n\frac{(b_jx_i/a_ix_j,qx_i/x_j;q)_{\infty}}
{(qx_i/a_ix_j,b_jx_i/x_j;q)_{\infty}},
\end{multline}
where $|q|<1$ and $|b_1\cdots b_nq^{1-n}/a_1\cdots a_n|<|z|<1$.

The special case of \eqref{1psi1gus}, in which $b_1 = \cdots = b_n = b$,
was previously obtained in \cite[Theorem~1.15]{M86}.
Further applications of this first multilateral ${}_{1}\psi_{1}$
summation appear in \cite{M86,M5,M6}.

Another $\mathrm A_n$ ${}_1\psi_1$ summation theorem was found in
\cite[Theorem~3.2]{MS}:
\begin{multline}\label{1psi1b}
\sum_{k_1,\dots,k_n=-\infty}^{\infty}\Bigg(
\prod_{1\le i<j\le n}\frac{x_iq^{k_i}-x_jq^{k_j}}{x_i-x_j}
\prod_{i,j=1}^n(b_jx_i/x_j;q)_{k_i}^{-1}
\prod_{i=1}^nx_i^{nk_i-|k|}\\\times
(a;q)_{|k|}(-1)^{(n-1)|k|}q^{-\binom{|k|}2+n\sum_{i=1}^n\binom{k_i}2}z^{|k|}\Bigg)\\
=\frac{(az,q/az,b_1\cdots b_nq^{1-n}/a;q)_\infty}
{(z,b_1\cdots b_nq^{1-n}/az,q/a;q)_\infty}
\prod_{i,j=1}^n\frac{(qx_i/x_j;q)_{\infty}}{(b_jx_i/x_j;q)_{\infty}},
\end{multline}
where $|q|<1$ and $|b_1\cdots b_nq^{1-n}/a|<|z|<1$.
(The specified region of convergence can be determined
by an analysis as carried out in \cite[Appendix~A]{Schl3}.)

Another $\mathrm A_n$ extension of Ramanujan's ${}_1\psi_1$ summation
is implicitly contained in \cite{Mac}.
Written out in explicit terms, it reads as \cite[last identity]{W2}
\begin{align}\label{mac11}
\sum_{k_1,\dots,k_n=-\infty}^{\infty}
\prod_{1\le i<j\le n}\Bigg(
\frac{x_iq^{k_i}-x_jq^{k_j}}{x_i-x_j}\,
\frac{(x_i/tx_j;q)_{k_i-k_j}}{(qtx_i/x_j;q)_{k_i-k_j}}
q^{-k_j}t^{k_i-k_j}\Bigg)\cdot\frac{(a;q)_{|k|}}{(b;q)_{|k|}}z^{|k|}&\notag\\
=\frac{(az,q/az,b/a,qt;q)_{\infty}}
{(z,b/az,q/a,b;q)_{\infty}}
\prod_{i=1}^{n-1}\frac{(qt^{i+1};q)_\infty}{(t^i;q)_\infty}
\prod_{i,j=1}^n\frac{(qx_i/x_j;q)_{\infty}}
{(qtx_i/x_j;q)_{\infty}}&,
\end{align}
where $|q|<1$, $|t|<1$ and $|b/a|<|z|<1$.

Taking coefficients of $z^N$ on both sides of \eqref{mac11}
while appealing to (the univariate version of) Ramanujan's
${}_1\psi_1$ summation, we obtain
the interesting identity
\begin{align}\label{mac11t}
\sum_{\substack{-\infty\le k_1,\dots,k_n\le\infty\\|k|=N}}
\prod_{1\le i<j\le n}\Bigg(
\frac{x_iq^{k_i}-x_jq^{k_j}}{x_i-x_j}\,
\frac{(x_i/tx_j;q)_{k_i-k_j}}{(qtx_i/x_j;q)_{k_i-k_j}}
q^{-k_j}t^{k_i-k_j}\Bigg)&\notag\\
=\frac{(qt;q)_{\infty}}{(q;q)_{\infty}}
\prod_{i=1}^{n-1}\frac{(qt^{i+1};q)_\infty}{(t^i;q)_\infty}
\prod_{i,j=1}^n\frac{(qx_i/x_j;q)_{\infty}}
{(qtx_i/x_j;q)_{\infty}}&,
\end{align}
subject to $|q|<1$. Observe that the right-hand side is independent of $N$.

Some additional (simpler) $\mathrm A_n$ ${}_1\psi_1$ summations are given
in \cite{GK1}, \cite[Section~2]{Schl1} and \cite[Section~6]{RS}.

Bailey's very-well-poised
${}_6\psi_6$ summation is (cf.\ \cite[Equation (II.33)]{GR})
\begin{align}
&{}_6\psi_6\!\left[\begin{matrix}qa^{\frac 12},-qa^{\frac 12},b,c,d,e\\
a^{\frac 12},-a^{\frac 12},aq/b,aq/c,aq/d,aq/e\end{matrix};q,q\right]\notag\\
&=\frac{(q,aq,q/a,aq/bc,aq/bd,aq/be,aq/cd,aq/ce,aq/de;q)_\infty}
{(aq/b,aq/c,aq/d,aq/e,q/b,q/c,q/d,q/e,a^2q/bcde;q)_\infty},
\end{align}
where $|q|<1$ and $|a^2q/bcde|<1$.
Several root system extensions of Bailey's ${}_6\psi_6$ summation
formula exist. Due to the fundamental importance of
the ${}_6\psi_6$ summation, we review several of these
summations:

We start with an $\mathrm A_n$ extension of the ${}_6\psi_6$ summation
\cite[Theorem~1.15]{G2}:
\begin{multline}\label{r66gl}
\sum_{k_1,\dots,k_n=-\infty}^{\infty}
\Bigg(\prod_{1\le i<j\le n}
\frac{x_iq^{k_i}-x_jq^{k_j}}{x_i-x_j}
\prod_{i=1}^n\frac{1-ax_iq^{k_i+|k|}}{1-ax_i}
\prod_{i,j=1}^n\frac{(b_jx_i/x_j;q)_{k_i}}
{(ax_iq/e_jx_j;q)_{k_i}}\\\times
\prod_{i=1}^n\frac{(e_ix_i;q)_{|k|}\,(cx_i;q)_{k_i}}
{(ax_iq/b_i;q)_{|k|}\,(ax_iq/d;q)_{k_i}}\cdot
\frac{(d;q)_{|k|}}{(aq/c;q)_{|k|}}\,
\left(\frac{a^{n+1}q}{BcdE}\right)^{|k|}\Bigg)\\
=\frac{(aq/Bc,a^nq/dE,aq/cd;q)_{\infty}}
{(a^{n+1}q/BcdE,aq/c,q/d;q)_{\infty}}
\prod_{i,j=1}^n\frac{(ax_iq/b_ie_jx_j,qx_i/x_j;q)_{\infty}}
{(qx_i/b_ix_j,ax_iq/e_jx_j;q)_{\infty}}\\\times
\prod_{i=1}^n\frac{(aq/ce_ix_i,ax_iq/b_id,ax_iq,q/ax_i;q)_{\infty}}
{(ax_iq/b_i,q/e_ix_i,q/cx_i,ax_iq/d;q)_{\infty}},
\end{multline}
where $B=b_1\dots b_n$ and $E=e_1\dots e_n$,
provided $|q|<1$ and $|a^{n+1}q/BcdE|<1$.

The multilateral identity above can also be
written in a more compact form. We then have the
$\mathrm A_n$ ${}_6\psi_6$ summation from \cite[Theorem~1.15]{G2}:
\begin{multline}\label{an1psi1cglN}
\sum_{\begin{smallmatrix}-\infty\le k_1,\dots,k_{n+1}\le\infty\\
k_1+\dots+k_{n+1}=0\end{smallmatrix}}
\prod_{1\le i<j\le n+1}\frac {x_iq^{k_i}-x_jq^{k_j}}{x_i-x_j}
\prod_{i,j=1}^{n+1}\frac{(a_jx_i/x_j;q)_{k_i}}{(b_jx_i/x_j;q)_{k_i}}\\
=\frac{(b_1\cdots b_{n+1}q^{-n},q/a_1\cdots a_{n+1};q)_{\infty}}
{(q,b_1\cdots b_{n+1}q^{-n}/a_1\cdots a_{n+1};q)_{\infty}}
\prod_{i,j=1}^{n+1}
\frac{(qx_i/x_j,b_jx_i/a_ix_j;q)_{\infty}}
{(b_jx_i/x_j,x_iq/a_ix_j;q)_{\infty}},
\end{multline}
provided $|q|<1$ and $|b_1\cdots b_{n+1}q^{-n}/a_1\cdots a_{n+1}|<1$.
It is not difficult to see that  \eqref{r66gl} and \eqref{an1psi1cglN}
are equivalent by a change of variables.

If in \eqref{an1psi1cglN} one replaces the parameters
$x_i$, $a_i$ and $b_i$, by $q^{x_i}$,  $q^{a_i}$ and  $q^{b_i}$, respectively,  
and formally lets $q\to 1$, one obtains the following
$\mathrm A_n$ ${}_5H_5$ summation from \cite[Theorem~1.13]{G2}
(given a direct proof there by functional equations
without appealing to a $q\to 1$ limit): 
\begin{multline}\label{an2h2cglN}
\sum_{\begin{smallmatrix}-\infty\le k_1,\dots,k_{n+1}\le\infty\\
k_1+\dots+k_{n+1}=0\end{smallmatrix}}
\prod_{1\le i<j\le n+1}\frac {x_i+k_i-x_j-k_j}{x_i-x_j}
\prod_{i,j=1}^{n+1}\frac{(a_j+x_i-x_j)_{k_i}}{(b_j+x_i-x_j)_{k_i}}\\
=
\frac{\Gamma\Big({-}n+\sum_{i=1}^{n+1}(b_i-a_i)\Big)}
{\Gamma\Big(1-\sum_{i=1}^{n+1}a_i\Big)\,
\Gamma\Big({-}n+\sum_{i=1}^{n+1}b_i\Big)}
\prod_{i,j=1}^{n+1}
\frac{\Gamma(b_j+x_i-x_j)\,\Gamma(1-a_i+x_i-x_j)}
{\Gamma(1+x_i-x_j)\,\Gamma(b_j-a_i+x_i-x_j)},
\end{multline}
provided $\mathrm{Re}\Big(\sum_{i=1}^{n+1}(b_i-a_i)\Big)>n$.

Another $\mathrm A_n$ very-well-poised
${}_6\psi_6$ summation was derived in \cite{Schl5}.
\begin{multline}
\sum_{k_1,\dots,k_n=-\infty}^{\infty}
\Bigg(\prod_{1\le i<j\le n}
\frac{x_iq^{k_i}-x_jq^{k_j}}{x_i-x_j}
\prod_{i,j=1}^n\frac{(c_jx_i/x_j;q)_{k_i}}{(ax_iq/e_jx_j;q)_{k_i}}\\\times
\prod_{i=1}^n\frac{(aq/bCx_i;q)_{|k|-k_i}\,(dE/a^{n-1}e_ix_i;q)_{|k|}\,
(bx_i;q)_{k_i}}
{(dE/a^nx_i;q)_{|k|-k_i}\,(ac_iq/bCx_i;q)_{|k|}\,(ax_iq/d;q)_{k_i}}\\\times
\frac{1-aq^{2|k|}}{1-a}\,\frac{(E/a^{n-1};q)_{|k|}}{(aq/C;q)_{|k|}}
\left(\frac{a^{n+1}q}{bCdE}\right)^{|k|}\Bigg)\\=
\frac{(aq,q/a,aq/bd;q)_\infty}{(aq/C,a^{n+1}q/bCdE,a^{n-1}q/E;q)_\infty}
\prod_{i,j=1}^n\frac{(qx_i/x_j,ax_iq/c_ie_jx_j;q)_\infty}
{(qx_i/c_ix_j,ax_iq/e_jx_j;q)_\infty}\\\times
\prod_{i=1}^n\frac{(a^nx_iq/dE,aq/be_ix_i,aq/bCx_i,ax_iq/c_id;q)_\infty}
{(a^{n-1}e_ix_iq/dE,q/bx_i,ax_iq/d,ac_iq/bCx_i;q)_\infty},
\end{multline}
where $C=c_1\cdots c_n$ and $E=e_1\cdots e_n$,
provided $|q|<1$ and $|a^{n+1}q/BcdE|<1$.

A $\mathrm C_n$ very-well-poised ${}_6\psi_6$ summation was
derived in \cite[Theorem~5.1]{G3}.
\begin{multline}\label{cr66gl}
\sum_{k_1,\dots,k_n=-\infty}^{\infty}
\Bigg(\prod_{1\le i<j\le n}\frac {x_iq^{k_i}-x_jq^{k_j}}{x_i-x_j}
\prod_{1\le i\le j\le n}\frac {1-ax_ix_jq^{k_i+k_j}}{1-ax_ix_j}\\\times
\prod_{i,j=1}^n\frac{(c_jx_i/x_j,e_jx_ix_j;q)_{k_i}}
{(ax_ix_jq/c_j,ax_iq/e_jx_j;q)_{k_i}}
\prod_{i=1}^n\frac{(bx_i,dx_i;q)_{k_i}}{(ax_iq/b,ax_iq/d;q)_{k_i}}\cdot
\left(\frac{a^{n+1}q}{bCdE}\right)^{|k|}\Bigg)\\
=\prod_{1\le i<j\le n}(ax_ix_jq/c_ic_j,aq/e_ie_jx_ix_j;q)_{\infty}
\prod_{1\le i\le j\le n}(ax_ix_jq,q/ax_ix_j;q)_{\infty}\\\times
\frac{(aq/bd;q)_{\infty}}{(a^{n+1}q/bCdE;q)_{\infty}}
\prod_{i,j=1}^n\frac{(ax_iq/c_ie_jx_j,qx_i/x_j;q)_{\infty}}
{(ax_iq/e_jx_j,q/e_jx_ix_j,ax_ix_jq/c_i,qx_i/c_ix_j;q)_{\infty}}\\\times
\prod_{i=1}^n\frac{(ax_iq/bc_x,aq/be_ix_i,ax_iq/c_id,aq/de_ix_i;q)_{\infty}}
{(ax_iq/b,q/bx_i,ax_iq/d,q/dx_i;q)_{\infty}},
\end{multline}
where $C=c_1\cdots c_n$ and $E=e_1\cdots e_n$,
provided $|q|<1$ and $|a^{n+1}q/bCdE|<1$.

Here is an  $\mathrm B_n^\vee$ (in Macdonald's~\cite{Mac72}
terminology for affine root systems; or labeled 
$\mathrm A_{2n-1}^{(2)}$ by Kac~\cite{K})
very-well-poised ${}_6\psi_6$ summation,
obtained in \cite[Theorem~6.1]{G3}. 
\begin{multline}\label{br66gl}
\sum_{\substack{-\infty\le k_1,\dots,k_n\le\infty\\|k|\,\equiv\, \sigma\;\,(\!\!\!\!\!\!\mod 2)}}
\Bigg(\prod_{1\le i<j\le n}\frac {x_iq^{k_i}-x_jq^{k_j}}{x_i-x_j}
\prod_{1\le i\le j\le n}\frac {1-ax_ix_jq^{k_i+k_j}}{1-ax_ix_j}\\\times
\prod_{i,j=1}^n\frac{(c_jx_i/x_j,e_jx_ix_j;q)_{k_i}}
{(ax_ix_jq/c_j,ax_iq/e_jx_j;q)_{k_i}}\cdot
\left(-\frac{a^n}{bCdE}\right)^{|k|}\Bigg)\\
=\prod_{1\le i<j\le n}(ax_ix_jq/c_ic_j,aq/e_ie_jx_ix_j;q)_{\infty}
\prod_{1\le i\le j\le n}(ax_ix_jq,q/ax_ix_j;q)_{\infty}\\\times
\frac{(-q;q)_{\infty}}{(-a^n/CE;q)_{\infty}}
\prod_{i,j=1}^n\frac{(ax_iq/c_ie_jx_j,qx_i/x_j;q)_{\infty}}
{(ax_iq/e_jx_j,q/e_jx_ix_j,ax_ix_jq/c_i,qx_i/c_ix_j;q)_{\infty}}\\\times
\prod_{i=1}^n\frac{(aqx_i^2/c_i^2,aq/e_i^2u_i^2;q^2)_{\infty}}
{(aqx_i^2,q/ax_i^2;q^2)_{\infty}},
\end{multline}
where $C=c_1\cdots c_n$ and $E=e_1\cdots e_n$, and where $\sigma=0,1$,
provided $|q|<1$ and $|a^n/CE|<1$.

As observed in \cite{SW}, the identity in \eqref{br66gl} is closely
connected to the $b=\sqrt{aq}$, $d=-\sqrt{aq}$ case of
the identity in \eqref{cr66gl}, where the sum
evaluates to twice the product on the right-hand side of \eqref{br66gl}
(the latter being independent of $\sigma$).

Another $\mathrm C_n$ very-well-poised ${}_6\psi_6$ summation was established
in \cite[Equation~(2.22)]{vD}.
\begin{multline}\label{cnvd66gl}
\sum_{k_1,\dots,k_n=-\infty}^{\infty}
\Bigg(\prod_{1\le i<j\le n}\frac {x_iq^{k_i}-x_jq^{k_j}}{x_i-x_j}
\prod_{1\le i\le j\le n}\frac {1-ax_ix_jq^{k_i+k_j}}{1-ax_ix_j}\\\times
\prod_{1\le i\le j\le n}\frac{(tax_ix_j;q)_{k_i+k_j}\,(tx_i/x_j;q)_{k_i-k_j}}
{(ax_ix_jq/t;q)_{k_i+k_j}\,(qx_i/tx_j;q)_{k_i-k_j}}\cdot
\left(\frac{t^2}q\right)^{\sum_{i=1}^n(i-1)k_i}\\\times
\prod_{i=1}^n\frac{(bx_i,cx_i,dx_i,ex_i;q)_{k_i}}
{(ax_iq/b,ax_iq/c,ax_iq/d,ax_iq/e;q)_{k_i}}\cdot
\left(\frac{t^{2-2n}a^2q}{bcde}\right)^{|k|}\Bigg)\\
=\prod_{i,j=1}^n\frac{(qx_i/x_j;q)_\infty}{(qx_i/tx_j;q)_\infty}
\frac{\prod_{1\le i\le j\le n}(ax_ix_jq,q/ax_ix_j;q)_{\infty}}
{\prod_{1\le i<j\le n}(ax_ix_jq/t,q/tax_ix_j;q)_{\infty}}
\prod_{i=1}^n\frac{(qt^{-i};q)_\infty}{(qt^{2-i-n}a^2/bcde;q)_\infty}\\\times
\prod_{i=1}^n\frac{(at^{1-i}q/bc,at^{1-i}q/bd,at^{1-i}q/be,
at^{1-i}q/cd,at^{1-i}q/ce,at^{1-i}q/de;q)_\infty}
{(q/bx_i,q/cx_i,q/dx_i,q/ex_i,
ax_iq/b,ax_iq/c,ax_iq/d,ax_iq/e;q)_\infty},
\end{multline}
provided $|q|<1$, $|a^2q^{2-n}/bcde|<1$ and $|t^{2-2n}a^2q/bcde|<1$.

The next  $\mathrm B_n^\vee$ (or $\mathrm A_{2n-1}^{(2)}$)
${}_6\psi_6$ summation (compare with \eqref{br66gl}) was derived in
\cite[Theorem~4.1]{SW}.
\begin{multline}\label{bn66gl}
\sum_{\substack{-\infty\le k_1,\dots,k_n\le\infty\\|k|\,\equiv\,\sigma\;\,(\!\!\!\!\!\!\mod 2)}}
\Bigg(\prod_{1\le i<j\le n}\frac {x_iq^{k_i}-x_jq^{k_j}}{x_i-x_j}
\prod_{1\le i\le j\le n}\frac {1-ax_ix_jq^{k_i+k_j}}{1-ax_ix_j}\\\times
\prod_{1\le i\le j\le n}\frac{(tax_ix_j;q)_{k_i+k_j}\,(tx_i/x_j;q)_{k_i-k_j}}
{(ax_ix_jq/t;q)_{k_i+k_j}\,(qx_i/tx_j;q)_{k_i-k_j}}\cdot
\left(\frac{t^2}q\right)^{\sum_{i=1}^n(i-1)k_i}\\\times
\prod_{i=1}^n\frac{(bx_i,cx_i;q)_{k_i}}
{(ax_iq/b,ax_iq/c;q)_{k_i}}\cdot
\left(-\frac{t^{2-2n}a}{bc}\right)^{|k|}\Bigg)\\
=\frac 12\prod_{i,j=1}^n\frac{(qx_i/x_j;q)_\infty}{(qx_i/tx_j;q)_\infty}
\frac{\prod_{1\le i\le j\le n}(ax_ix_jq,q/ax_ix_j;q)_{\infty}}
{\prod_{1\le i<j\le n}(ax_ix_jq/t,q/tax_ix_j;q)_{\infty}}
\prod_{i=1}^n\frac{(qt^{-i};q)_\infty}{(-t^{2-i-n}a/bc;q)_\infty}\\\times
\prod_{i=1}^n\frac{(at^{1-i}q/bc,-t^{1-i};q)_\infty\,
(at^{2-2i}q/b,at^{2-2i}q/c;q^2)_\infty}
{(q/bx_i,q/cx_i,ax_iq/b,ax_iq/c;q)_\infty\,
(q/ax_i^2,aqx_i^2;q^2)_\infty},
\end{multline}
 where $\sigma=0,1$, provided $|q|<1$, $|aq^{1-n}/bc|<1$ and $|t^{2-2n}a/bc|<1$.

Similar to the relation between \eqref{br66gl} and \eqref{cr66gl},
the identity in \eqref{bn66gl} is closely connected to
the $d=\sqrt{aq}$, $e=-\sqrt{aq}$ case of the identity in \eqref{cnvd66gl},
where the sum evaluates to twice the product on the right-hand
side of \eqref{bn66gl}
(the latter being independent of $\sigma$).

Two other (simpler) $\mathrm C_n$ ${}_6\psi_6$ summations are given in
\cite[Theorem~3.4]{Schl1}.

Multivariate analogues of Bailey's ${}_6\psi_6$ summation for
\textit{exceptional} root systems were derived in
\cite{G4} (summation for $\mathrm G_2$),
\cite{I1} (summation for $\mathrm F_4$; see also \cite{I2}),
and \cite{IT} (further summations for $\mathrm G_2$).

\subsection{Watson transformations}

The Watson transformation (cf.\ \cite[Equation (III.18)]{GR})
\begin{align}
&{}_8\phi_7\!\left[\begin{matrix}a,\,qa^{\frac 12},-qa^{\frac 12},b,c,d,e,q^{-N}\\
a^{\frac 12},-a^{\frac 12},aq/b,aq/c,aq/d,aq/e,aq^{N+1}
\end{matrix};q,q\right]\notag\\
&=\frac{(aq,aq/de;q)_N}
{(aq/d,aq/e;q)_N}\,{}_4\phi_3\!\left[\begin{matrix}aq/bc,d,e,q^{-N}\\
aq/b,aq/c,deq^{-N}/a
\end{matrix};q,q\right],
\end{align}
where $|a^2q^{N+2}/bcde|<1$,
is very useful. For instance,
it can be used for a quick proof of the Rogers--Ramanujan identities,
see \cite[Section~2.7]{GR}.

A number of Watson transformations for basic hypergeometric series
associated with root systems exist. We only reproduce a few of them here.

In \cite[Theorem~2.24]{G1}, Gustafson applied the representation theory
of $\mathrm U(n)$ to derive the first multivariable generalization
of Whipple's classical transformation of an ordinary ($q = 1$)
terminating well-poised ${}_7F_6(1)$ into a terminating balanced
${}_4F_3(1)$. Its $q$-analogue, the first multivariable Watson
transformation, was obtained in \cite[Theorems~6.1 and 6.4]{M7}
and \cite[Theorems~6.1 and 6.4]{M89}
by a direct, elementary proof utilizing $q$-difference equations and
induction. Further details and applications are given in
\cite{M7,M89,M94}. A more symmetrical $\mathrm A_n$ Watson
transformation was derived in \cite[Theorem~2.1]{M00},
by means of the summation theorems
and analysis from \cite{M6}, where \cite{M6} provides an
$\mathrm A_n$ generalization of much of the analysis in chapters
one and two of Gasper and Rahman\rq{}s book \cite{GR}.

The following $\mathrm A_n$ Watson transformation was derived in
\cite[Theorem~A.3]{MN}.
\begin{multline}\label{anwatsongl}
\sum_{k_1,\dots,k_n=0}^{N_1,\dots,N_n}\Bigg(
\prod_{1\le i<j\le n}\frac{x_iq^{k_i}-x_jq^{k_j}}{x_i-x_j}
\prod_{i,j=1}^n\frac{(q^{-N_j}x_i/x_j;q)_{k_i}}
{(qx_i/x_j;q)_{k_i}}\prod_{i=1}^n\frac{1-ax_iq^{k_i+|k|}}
{1-ax_i}\,\frac{(ax_i;q)_{|k|}}{(ax_iq^{1+N_i};q)_{|k|}}\\\times
\prod_{i=1}^n\frac{(bx_i,cx_i;q)_{k_i}}
{(ax_iq/d,ax_iq/e;q)_{k_i}}\cdot
\frac{(d,e;q)_{|k|}}{(aq/b,aq/c;q)_{|k|}}
\left(\frac{a^2q^{|N|+2}}{bcde}\right)^{|k|}\Bigg)\\=
\frac{(aq/ce;q)_{|N|}}{(aq/c;q)_{|N|}}
\prod_{i=1}^n\frac{(ax_iq;q)_{N_i}}{(ax_iq/e;q)_{N_i}}\\\times
\sum_{k_1,\dots,k_n=0}^{N_1,\dots,N_n}\Bigg(
\prod_{1\le i<j\le n}\frac{x_iq^{k_i}-x_jq^{k_j}}{x_i-x_j}
\prod_{i,j=1}^n\frac{(q^{-N_j}x_i/x_j;q)_{k_i}}
{(qx_i/x_j;q)_{k_i}}\\\times
\prod_{i=1}^n\frac{(cx_i;q)_{k_i}}{(ax_iq/d;q)_{k_i}}\cdot
\frac{(aq/bd,e;q)_{|k|}}{(aq/b,ceq^{-|N|}/a;q)_{|k|}}q^{|k|}\Bigg).
\end{multline}
For a very similar but different $\mathrm A_n$ Watson transformation, see the
$f_i=q^{-N_i}$, $i=1,\dots,n$, case of \eqref{anntwatsongl}.

The following $\mathrm C_n\leftrightarrow\mathrm A_{n-1}$ 
Watson transformation was first derived in
\cite[Theorem~6.6]{ML}.
\begin{multline}\label{cnwatsongl}
\sum_{k_1,\dots,k_n=0}^{N_1,\dots,N_n}\Bigg(
\prod_{1\le i<j\le n}\frac{x_iq^{k_i}-x_jq^{k_j}}{x_i-x_j}
\prod_{1\le i\le j\le n}\frac{1-x_ix_jq^{k_i+k_j}}{1-x_ix_j}
\prod_{i,j=1}^n\frac{(q^{-N_j}x_i/x_j,x_ix_j;q)_{k_i}}
{(qx_i/x_j,q^{1+N_i}x_ix_j;q)_{k_i}}\\\times
\prod_{i=1}^n\frac{(bx_i,cx_i,dx_i,ex_i;q)_{k_i}}
{(qx_i/b,qx_i/c,qx_i/d,qx_i/e;q)_{k_i}}\cdot
\left(\frac{q^{|N|+2}}{bcde}\right)^{|k|}\Bigg)\\
=(q/bc;q)_{|N|}\prod_{i=1}^n\frac 1{(qx_i/b,qx_i/c;q)_{N_i}}
\prod_{i,j=1}^n(qx_ix_j;q)_{N_i}\prod_{1\le i<j\le n}
\frac 1{(qx_ix_j;q)_{N_i+N_j}}\\\times
\sum_{k_1,\dots,k_n=0}^{N_1,\dots,N_n}\Bigg(
\prod_{1\le i<j\le n}\frac{x_iq^{k_i}-x_jq^{k_j}}{x_i-x_j}
\prod_{i,j=1}^n\frac{(q^{-N_j}x_i/x_j;q)_{k_i}}
{(qx_i/x_j;q)_{k_i}}\\\times
\prod_{i=1}^n\frac{(bx_i,cx_i;q)_{k_i}}
{(qx_i/d,qx_i/e;q)_{k_i}}\cdot
\frac{(q/de;q)_{|k|}}{(bcq^{-|N|};q)_{|k|}}q^{|k|}\Bigg).
\end{multline}
This identity was utilized in \cite{Ba} and in
\cite{BW} to obtain identities for characters of
affine Lie algebras.

Several other Watson transformations are given in \cite{B,BS,BS2,C}.
One of them is the following (cf.\ \cite[Theorem~4.10]{BS}):
\begin{multline}\label{dnwat1e}
\sum_{k_1,\dots,k_n=0}^{N_1,\dots,N_n}\Bigg(
\prod _{1\le i<j\le n}\frac{x_iq^{k_i}-x_jq^{k_j}}
{x_i-x_j}\,\frac{(ax_ix_jq/c;q)_{k_i+k_j}}{(ex_ix_j;q)_{k_i+k_j}}
\prod_{i=1}^{n}\frac{1-ax_iq^{k_i+|k|}}{1-ax_i}\,
\frac{(aq/ex_i;q)_{|k|-k_i}}{(c/x_i;q)_{|k|-k_i}}\\\times
\prod_{i=1}^{n}\frac {(ax_i,c/x_i;q)_{|k|}\,(bx_i;q)_{k_i}}
{(ax_iq^{1+N_i},aq^{1-N_i}/ex_i;q)_{|k|}\,(ax_iq/d;q)_{k_i}}\\\times
\prod_{i,j=1}^{n}\frac{(q^{-N_j}x_i/x_j,ex_ix_jq^{N_j};q)_{k_i}}
{(qx_i/x_j,ax_ix_jq/c;q)_{k_i}}\cdot
\frac{(d;q)_{|k|}}{(aq/b;q)_{|k|}}
\left({\frac {q^2a^2}{bcde}}\right)^{|k|}\Bigg)\\
=d^{-|N|}\prod_{i=1}^{n}\frac{(ax_iq,dex_i/a;q)_{N_i}}
{(ex_i/a,ax_iq/d;q)_{N_i}}\\\times
\sum_{k_1,\dots,k_n=0}^{N_1,\dots,N_n}\Bigg(
\prod_{1\le i<j\le n}\frac{x_iq^{k_i}-x_jq^{k_j}}
{x_i-x_j}\,\frac{(ax_ix_jq/c;q)_{k_i+k_j}}{(ex_ix_j;q)_{k_i+k_j}}
\prod_{i=1}^{n}\frac{(ax_iq/bc;q)_{k_i}}{(dex_i/a;q)_{k_i}}\\\times
\prod_{i,j=1}^{n}\frac{(q^{-N_j}x_i/x_j,ex_ix_jq^{N_j};q)_{k_i}}
{(qx_i/x_j,ax_ix_jq/c;q)_{k_i}}\cdot
\frac{(d;q)_{|k|}}{(aq/b;q)_{|k|}}q^{|k|}\Bigg).
\end{multline}

This multivariate Watson transformation cannot be simplified to
any multivariate Jackson summation as a special case.

\subsection{Dimension changing transformations}

Heine's $q$-analogue of the classical
Euler transformation of ${}_2F_1$ series is (cf.\ \cite[Equation (III.3)]{GR})
\begin{equation}\label{qeuleru}
{}_2\phi_1\left[\begin{matrix}a,b\\c\end{matrix};q,z\right]=
\frac{(abz/c;q)_\infty}{(z;q)_\infty}\,
{}_2\phi_1\left[\begin{matrix}c/a,c/b\\c\end{matrix};
q,\frac{abz}c\right],
\end{equation}
valid for $|q|<1$, $|z|<1$ and $|abz/c|<1$.
The following result, which was first derived by Kajihara \cite{K1},
connects $\mathrm A_n$ and $\mathrm A_m$ basic hypergeometric series
and reduces, for $n=m=1$, to the $q$-Euler transformation.
\begin{align}\label{kaji1id}\notag
&\sum_{k_1,\dots,k_n\ge0}
\prod_{1\le i<j\le n}\frac{x_iq^{k_i}-x_jq^{k_j}}{x_i-x_j}
\prod_{1\le i,j\le n}\frac{(a_jx_i/x_j;q)_{k_i}}{(qx_i/x_j;q)_{k_i}}
\prod_{\substack{1\le i\le n\\1\le l\le m}}
\frac{(b_lx_iy_l;q)_{k_i}}{(cx_iy_l;q)_{k_i}}\cdot z^{|k|}\\&
=\frac{(ABz/c^m;q)_\infty}{(z;q)_\infty}\notag\\&\;\;\times
\sum_{\ka_1,\dots,\ka_m\ge0}
\prod_{1\le j<l\le m}\frac{y_jq^{\ka_j}-y_lq^{\ka_l}}{y_j-y_l}
\prod_{1\le j,l\le m}\frac{(cy_j/b_ly_l;q)_{\ka_j}}{(qy_j/y_l;q)_{\ka_j}}
\prod_{\substack{1\le i\le n\\1\le l\le m}}
\frac{(cx_iy_l/a_i;q)_{\ka_l}}{(cx_iy_l;q)_{\ka_l}}\cdot
\left(\frac
{ABz}{c^m}\right)^{|\ka|},
\end{align}
where $A=a_1\cdots a_n$, and $B=b_1\cdots b_m$, provided
$|q|<1$, $|z|<1$ and $|ABz/c^m|<1$.

Now let
\begin{align}
&\Phi_N^{n,m}\left(\left.\begin{matrix}\{a_i\}_n\\\{x_i\}_n\end{matrix}\right|
\begin{matrix}\{b_ly_l\}_m\\\{cy_l\}_m\end{matrix}\right)\notag\\
&:=
\sum_{\substack{k_1,\dots,k_n\ge0\\|k|=N}}
\prod_{1\le i<j\le n}\frac{x_iq^{k_i}-x_jq^{k_j}}{x_i-x_j}
\prod_{1\le i,j\le n}\frac{(a_jx_i/x_j;q)_{k_i}}{(qx_i/x_j;q)_{k_i}}
\prod_{\substack{1\le i\le n\\1\le l\le m}}
\frac{(b_lx_iy_l;q)_{k_i}}{(cx_iy_l;q)_{k_i}}.
\end{align}

The transformation in \eqref{kaji1id} was used to derive the following identity
\cite[Theorem 3.1]{K2} (which we state here in corrected form):
\begin{multline}
\sum_{K=0}^N
\Phi_{K}^{n_2,m_2}\left(\left.\begin{matrix}\{f/e_t\}_{n_2}\\
\{v_t\}_{n_2}\end{matrix}\right|
\begin{matrix}\{fw_r/d_r\}_{m_2}\\\{fw_r\}_{m_2}\end{matrix}\right)
\Phi_{N-K}^{n_1,m_1}\left(\left.\begin{matrix}\{a_i\}_{n_1}\\
\{x_i\}_{n_1}\end{matrix}\right|
\begin{matrix}\{b_ly_l\}_{m_1}\\\{cy_l\}_{m_1}\end{matrix}\right)
\left(\frac{d_1\cdots d_{m_2}e_1\cdots e_{n_2}}{f^{n_2}}\right)^{K}\\=
\sum_{L=0}^N
\Phi_L^{m_1,n_1}\left(\left.\begin{matrix}\{c/b_l\}_{m_1}\\
\{y_l\}_{m_1}\end{matrix}\right|
\begin{matrix}\{cx_i/a_i\}_{n_1}\\\{cx_i\}_{n_1}\end{matrix}\right)
\Phi_{N-L}^{m_2,n_2}\left(\left.\begin{matrix}\{d_r\}_{m_2}\\
\{w_r\}_{m_2}\end{matrix}\right|
\begin{matrix}\{e_tv_t\}_{n_2}\\\{fv_t\}_{n_2}\end{matrix}\right)
\left(\frac{a_1\cdots a_{n_1}b_1\cdots b_{m_1}}{c^{m_1}}\right)^L,
\end{multline}
where $a_1\cdots a_{n_1}b_1\cdots b_{m_1}/c^{m_1}=
d_1\cdots d_{m_2}e_1\cdots e_{n_2}/f^{n_2}$.
This identity can be viewed as a multivariate
extension of the Sears transformation \cite[Equation (III.16)]{GR})
(to which it reduces for $n=m=1$ after some elementary manipulations).
A transformation similar to \eqref{kaji1id} but connecting
$\mathrm C_n$ and $\mathrm C_m$ basic hypergeometric series
has been given in \cite{KMN}.

Several other transformations connecting sums of different dimension exist.
For instance, in \cite{R4} the following reduction formula for a multilateral
Karlsson--Minton type basic hypergeometric series associated with the
root system $\mathrm A_n$ was derived. (A basic hypergeometric series
is said to be of {\it Karlsson--Minton type} if the quotient
of corresponding upper and lower parameters is a nonnegative integer
power of $q$.)
\begin{multline}\label{rosanid}
\sum_{\substack{k_1,\dots,k_n=-\infty\\k_1+\dots+k_n=0}}^\infty
\prod_{1\leq i<j\leq n}\frac{x_iq^{k_i}-x_jq^{k_j}}{x_i-x_j}
\prod_{\substack{1\leq i\leq n\\1\leq j\leq p}}\frac{(x_iy_jq^{m_j};q)_{k_i}}
{(x_iy_j;q)_{k_i}}\prod_{i,j=1}^n\frac{(x_ia_j;q)_{k_i}}
{(x_ib_j;q)_{k_i}}\\
=\frac{(q^{1-|m|}/AX,
q^{1-n}BX;q)_\infty}
{(q,q^{1-|m|-n}B/A;q)_\infty}
\prod_{i,j=1}^n\frac{(b_i/a_j,qx_i/x_j;q)_\infty}
{(q/x_ia_j,x_ib_j;q)_\infty}
\prod_{\substack{1\leq i\leq n\\1\leq j\leq p}}\frac{(q^{-m_j}b_i/y_j;q)_{m_j}}
{(q^{1-m_j}/x_iy_j;q)_{m_j}}\\
\times
\sum_{\ka_1,\dots,\ka_p=0}^{m_1,\dots,m_p}q^{|\ka|}
\frac{(q^{n}/BX;q)_{|\ka|}}
{(q^{1-|m|}/AX;q)_{|\ka|}}
\prod_{1\leq i<j\leq n}\frac{y_iq^{\ka_i}-y_jq^{\ka_j}}{y_i-y_j}
\prod_{\substack{1\leq i\leq n\\1\leq j\leq p}}\frac{(y_j/a_i;q)_{\ka_j}}
{(qy_j/b_i;q)_{\ka_j}}\prod_{i,j=1}^p\frac{(q^{-m_i}y_j/y_i;q)_{\ka_j}}
{(qy_j/y_i;q)_{\ka_j}},
\end{multline}
where $A=a_1\cdots a_n$, $B=b_1\cdots b_n$, $X=x_1\cdots x_n$, provided
$|q|<1$ and $|q^{1-|m|-n}B/A|<1$.

Similarly, in \cite{R3} the following reduction formula for a multilateral
Karlsson--Minton type basic hypergeometric series associated with the
root system $\mathrm C_n$ was derived.
\begin{multline}\label{roscnid}
\sum_{k_1,\dots,k_n=-\infty}^\infty\Bigg(
\prod_{1\leq i<j\leq n}\frac{x_iq^{k_i}-x_jq^{k_j}}{x_i-x_j}
\prod_{1\leq i\leq j\leq n}
\frac{1-x_ix_jq^{k_i+k_j}}{1-x_ix_j}\\\times
\prod_{\substack{1\leq i\leq n\\1\leq j\leq p}}
\frac{(x_iy_jq^{m_j},qx_i/y_j;q)_{k_i}}
{(x_iy_j,q^{1-m_j}x_i/y_j;q)_{k_i}}
\prod_{\substack {1\leq i\leq n\\1\leq j\leq 2n+2}}
\frac{(x_ia_j;q)_{k_i}}{(qx_i/a_j;q)_{k_i}}\cdot
\left(\frac{q^{1-|m|}}A\right)^{|k|}\Bigg)\\
=\frac{\prod_{1\leq i\leq j\leq n}(qx_ix_j,q/x_ix_j;q)_\infty
\prod_{i,j=1}^n(qx_i/x_j;q)_\infty}
{\prod_{ {1\leq i\leq n,\,1\leq j\leq 2n+2}}(qx_i/a_j,q/x_ia_j;q)_\infty}
\frac{\prod_{1\leq i<j\leq 2n+2}(q/a_ia_j;q)_\infty}
{(q/A;q)_\infty}\\\times
\frac{\prod_{{1\leq i\leq 2n+2,\,1\leq j\leq p}}
(y_ja_i;q)_{m_j}}
{\prod_{{1\leq i\leq n,\,1\leq j\leq p}}
(y_jx_i,y_j/x_i;q)_{m_j}}
\frac{\prod_{1\leq i<j\leq p}(y_iy_j;q)_{m_i+m_j}}
{\prod_{i,j=1}^p(y_iy_j;q)_{m_i}}\frac 1{(A;q)_{|m|}}\\\times
\sum_{\ka_1,\dots,\ka_p=0}^{m_1,\dots,m_p}\Bigg(
\prod_{1\leq i<j\leq n}\frac{y_iq^{\ka_i}-y_jq^{\ka_j}}{y_i-y_j}
\prod_{1\leq i\leq j\leq n}
\frac{1-y_iy_jq^{\ka_i+\ka_j-1}}{1-y_iy_jq^{-1}}\\\times
\prod_{\substack {1\leq i\leq 2n+2\\1\leq j\leq p}}\frac{(y_j/a_i;q)_{\ka_j}}
{(y_ja_i;q)_{\ka_j}}\prod_{i,j=1}^p\frac{(q^{-1}y_iy_j,q^{-m_j}y_i/y_j;q)_{\ka_i}}
{(qy_i/y_j,q^{m_j}y_iy_j;q)_{\ka_i}}\cdot
\left(Aq^{|m|}\right)^{|\ka|}\Bigg),
\end{multline}
where $A=a_1\cdots a_{2n+2}$, provided $|q|<1$ and $|q^{1-|m|}/A|<1$.
A substantially more general transformation (involving four-fold
multiple sums) was given by Masuda~\cite[Theorem~3]{Ma}.

Both \eqref{rosanid} and \eqref{roscnid} have many interesting
consequences. In particular, these transformations form bridges between the
one-variable and the multivariable theory and can be used to prove various
summations and transformations for $\mathrm A_n$ and $\mathrm C_n$
basic hypergeometric series.
For details, we refer the reader to Rosengren's papers \cite{R3,R4}.

Other dimension changing transformations have been given (or conjectured)
in \cite{B1,GK,Kr,Ra1,R9}.

\subsection{Multiterm transformations}

Bailey's nonterminating balanced very-well-poised ${}_{10}\phi_9$ transformation
is (cf.\ \cite[Equation (III.39)]{GR})
\begin{align}\label{4t109}
&{}_{10}\phi_9\!\left[\begin{matrix}a,\,qa^{\frac 12},-qa^{\frac 12},
b,c,d,e,f,g,h\\a^{\frac 12},-a^{\frac 12},
aq/b,aq/c,aq/d,aq/e,aq/f,aq/g,aq/h\end{matrix}\,;q,q\right]\notag\\
&+\frac{(aq,b/a,c,d,e,f,g,h,bq/c,bq/d,bq/e,bq/f,bq/g,bq/h;q)_{\infty}}
{(b^2q/a,a/b,aq/c,aq/d,aq/e,aq/f,aq/g,aq/h,
bc/a,bd/a,be/a,bf/a,bg/a,bh/a;q)_{\infty}}
\notag\\&\quad\;\times
{}_{10}\phi_9\!\left[\begin{matrix}b^2/a,\,qba^{-\frac 12},-qba^{-\frac 12},
b,bc/a,bd/a,be/a,bf/a,bg/a,bh/a\\
ba^{-\frac 12},-ba^{-\frac 12},
bq/a,bq/c,bq/d,bq/e,bq/f,bq/g,bq/h\end{matrix}\,;q,q\right]\notag\\
&=\frac{(aq,b/a,\lambda q/f,\lambda q/g,\lambda q/h,
bf/\lambda,bg/\lambda,bh/\lambda;q)_{\infty}}
{(\lambda q,b/\lambda,aq/f,aq/g,aq/h,bf/a,bg/a,bh/a;q)_{\infty}}\notag\\
&\quad\;\times {}_{10}\phi_9\!\left[\begin{matrix}\lambda,\,q\lambda^{\frac 12},
-q\lambda^{\frac 12},b,\lambda c/a,\lambda d/a,\lambda e/a,f,g,h\\
\lambda^{\frac 12},-\lambda^{\frac 12},\lambda q/b,
aq/c,aq/d,aq/e,\lambda q/f,\lambda q/g,\lambda q/h\end{matrix}
\,;q,q\right]\notag\\
&\quad+\frac{(aq,b/a,f,g,h,bq/f,bq/g,bq/h,\lambda c/a,\lambda d/a,\lambda e/a,
abq/\lambda c,abq/\lambda d,abq/\lambda e;q)_{\infty}}
{(b^2q/\lambda,\lambda/b,aq/c,aq/d,aq/e,aq/f,aq/g,aq/h,
bc/a,bd/a,be/a,bf/a,bg/a,bh/a;q)_{\infty}}\notag\\
&\qquad\;\times
{}_{10}\phi_9\!\left[\begin{matrix}b^2/\lambda,\,qb\lambda^{-\frac 12},
-qb\lambda^{-\frac 12},b,bc/a,bd/a,be/a,bf/\lambda,bg/\lambda,bh/\lambda\\
b\lambda^{-\frac 12},-b\lambda^{-\frac 12},bq/\lambda,
abq/c\lambda,abq/d\lambda,abq/e\lambda,bq/f,bq/g,bq/h\end{matrix}\,;q,q\right],
\end{align}
where $\lambda=a^2q/cde$, $a^3q^2=bcdefgh$ and $|q|<1$.
This identity  is at the top of the classical
hierarchy of identities for basic hypergeometric series.

The following identity from \cite[Corollar 4.1]{RS} (which was derived
by determinant evaluations, following a method first used by
Gustafson and Krattenthaler~\cite{GK1,GK2} to derive $\mathrm A_n$ extensions
of Heine's ${}_2\phi_1$ transformations,
and subsequently used in a systematic manner in \cite{Schl1}
and \cite{RS})
concerns a $\mathrm C_n$
extension of Bailey's four-term transformation in \eqref{4t109},
where both sides of the identity involve
$2^n$ nonterminating $\mathrm C_n$ basic hypergeometric series.
Let $a^3q^{3-n}=bc_id_ie_ix_ifgh$ and
$\lambda=a^2q/c_id_ie_ix_i$ for $i=1,\dots,n$.
Then there holds
\begin{multline}\label{cnnt109gl}
\sum_{S\subseteq\{1,2,\dots,n\}}\Bigg[
\left(\frac ba\right)^{\binom{n-|S|}2}
\prod_{i\notin S}
\frac{(ax_i^2q,c_ix_i,d_ix_i,e_ix_i;q)_{\infty}}
{(ax_i/b,ax_iq/c_i,ax_iq/d_i,ax_iq/e_i;q)_{\infty}}\\\times
\prod_{i\notin S}
\frac{(fx_i,gx_i,hx_i,b/ax_i,bq/c_i,bq/d_i,bq/e_i,bq/f,bq/g,bq/h;q)_{\infty}}
{(ax_iq/f,ax_iq/g,ax_iq/h,b^2q/a,bc_i/a,bd_i/a,be_i/a,
bf/a,bg/a,bh/a;q)_{\infty}}\\\times
\sum_{k_1,\dots,k_n=0}^{\infty}\Bigg(
\underset{i,j\in S}{\prod_{1\le i<j\le n}}
\frac{(x_iq^{k_i}-x_jq^{k_j})(1-ax_ix_jq^{k_i+k_j})}
{(x_i-x_j)(1-ax_ix_j)}
\prod_{i\in S}\frac{1-ax_i^2q^{2k_i}}{1-ax_i^2}\\\times
\underset{i,j\notin S}{\prod_{1\le i<j\le n}}
\frac{(q^{k_i}-q^{k_j})(1-b^2q^{k_i+k_j}/a)}
{(x_i-x_j)(1-ax_ix_j)}
\prod_{i\notin S}\frac{1-b^2q^{2k_x}/a}{1-b^2/a}\\\times
\prod_{i\in S, j\notin S}
\frac{(x_iq^{k_i}-bq^{y_j}/a)(1-bx_jq^{k_i+k_j})}
{(x_i-x_j)(1-ax_ix_j)}\\\times
\prod_{i\in S}\frac{(ax_i^2,bx_i,c_ix_i,d_ix_i,e_ix_i,fx_i,gx_i,hx_i;q)_{k_i}}
{(q,ax_iq/b,ax_iq/c_i,ax_iq/d_i,ax_iq/e_i,
ax_iq/f,ax_iq/g,ax_iq/h;q)_{k_i}}\\\times
\prod_{i\notin S}\frac{(b^2/a,bx_i,bc_i/a,bd_i/a,be_i/a,
bf/a,bg/a,bh/a;q)_{k_i}}
{(q,bq/ax_i,bq/c_i,bq/d_i,bq/e_i,bq/f,bq/g,bq/h;q)_{k_i}}\,
\cdot q^{|k|}\Bigg)\Bigg]\\
=\prod_{i=1}^n\frac{(ax_i^2q,b/ax_i,
\lambda x_iq/f,\lambda x_iq/g,\lambda x_iq/h,
bfq^{i-1}/\lambda,bgq^{i-1}/\lambda,bhq^{i-1}/\lambda;q)_{\infty}}
{(\lambda x_i^2q,b/\lambda x_i,ax_iq/f,ax_iq/g,ax_iq/h,
bfq^{i-1}/a,bgq^{i-1}/a,bhq^{i-1}/a;q)_{\infty}}\\\times
\prod_{1\le i<j\le n}\frac{1-\lambda x_ix_j}{1-ax_ix_j}\;
\sum_{S\subseteq\{1,2,\dots,n\}}\Bigg[
\left(\frac b{\lambda}\right)^{\binom{n-|S|}2}\\\times
\prod_{i\notin S}
\frac{(\lambda x_i^2q,\lambda c_ix_i/a,\lambda d_ix_i/a,\lambda e_ix_i/a,
fx_i,gx_i,hx_i;q)_{\infty}}
{(\lambda x_i/b,ax_iq/c_i,ax_iq/d_i,ax_iq/e_i,
\lambda x_iq/f,\lambda x_iq/g,\lambda x_iq/h;q)_{\infty}}\\\times
\prod_{i\notin S}
\frac{(b/\lambda x_i,abq/c_i\lambda,abq/d_i\lambda,abq/e_i\lambda,
bq/f,bq/g,bq/h;q)_{\infty}}
{(b^2q/\lambda,bc_i/a,bd_i/a,be_i/a,
bf/\lambda,bg/\lambda,bh/\lambda;q)_{\infty}}\\\times
\sum_{k_1,\dots,k_n=0}^{\infty}\Bigg(
\underset{i,j\in S}{\prod_{1\le i<j\le n}}
\frac{(x_iq^{k_i}-x_jq^{k_j})(1-\lambda x_ix_jq^{k_i+k_j})}
{(x_i-x_j)(1-\lambda x_ix_j)}
\prod_{i\in S}\frac{1-\lambda x_i^2q^{2k_i}}{1-\lambda x_i^2}\\\times
\underset{i,j\notin S}{\prod_{1\le i<j\le n}}
\frac{(q^{k_i}-q^{k_j})(1-b^2q^{k_i+k_j}/\lambda)}
{(x_i-x_j)(1-\lambda x_ix_j)}
\prod_{i\notin S}\frac{1-b^2q^{2k_i}/\lambda}{1-b^2/\lambda}\\\times
\prod_{i\in S, j\notin S}
\frac{(x_iq^{k_i}-bq^{k_j}/\lambda)(1-bx_iq^{k_i+k_j})}
{(x_i-x_j)(1-\lambda x_ix_j)}\\\times
\prod_{i\in S}\frac{(\lambda x_i^2,bx_i,\lambda c_ix_i/a,
\lambda d_ix_i/a,\lambda e_ix_i/a,fx_i,gx_i,hx_i;q)_{k_i}}
{(q,\lambda x_iq/b,ax_iq/c_i,ax_iq/d_i,ax_iq/e_i,
\lambda x_iq/f,\lambda x_iq/g,\lambda x_iq/h;q)_{k_i}}\\\times
\prod_{i\notin S}\frac{(b^2/\lambda,bx_i,bc_i/a,bd_i/a,be_i/a,
bf/\lambda,bg/\lambda,bh/\lambda;q)_{k_i}}
{(q,bq/\lambda x_i,abq/c_i\lambda,abq/d_i\lambda,abq/e_i\lambda,
bq/f,bq/g,bq/h;q)_{k_i}}\,
\cdot q^{|k|}\Bigg)\Bigg],
\end{multline}
where $|q|<1$.

The next identity from \cite{MN2} concerns an
$\mathrm A_n$ extension of the nonterminating Watson transformation
(cf.\ \cite[Equation (III.36)]{GR}),
\begin{align}\label{3t87}
&{}_8\phi_7\!\left[\begin{matrix}a,\,qa^{\frac 12},-qa^{\frac 12},
b,c,d,e,f\\a^{\frac 12},-a^{\frac 12},
aq/b,aq/c,aq/d,aq/e,aq/f\end{matrix}\,;q,\frac{a^2q^2}{bcdef}\right]\notag\\
&=\frac{(aq,aq/de,aq/df,aq/ef;q)_{\infty}}
{(aq/d,aq/e,aq/f,aq/def;q)_{\infty}}\,
{}_4\phi_3\!\left[\begin{matrix}aq/bc,d,e,f\\
aq/b,aq/c,def/a\end{matrix}
\,;q,q\right]\notag\\
&\quad+\frac{(aq,aq/bc,d,e,f,a^2q^2/bdef,a^2q^2/cdef;q)_{\infty}}
{(aq/b,aq/c,aq/d,aq/e,aq/f,a^2q^2/bcdef,def/aq;q)_{\infty}}\notag\\
&\qquad\;\times
{}_4\phi_3\!\left[\begin{matrix}aq/de,aq/df,aq/ef,a^2q^2/bcdefa\\
a^2q^2/bdef,a^2q^2/cdef,aq^2/def\end{matrix}\,;q,q\right],
\end{align}
where $|q|<1$ and $|a^2q^2/bcdef|<1$.
In the multivariate case this
is a transformation of a nonterminating very-well-poised $\mathrm A_n$
basic hypergeometric series into $n+1$ multiples of
nonterminating balanced ${}\mathrm A_n$ basic hypergeometric series.
It is interesting to point out that although the $\mathrm A_n$ ${}_8\phi_7$
series on the left-hand side of \eqref{anntwatsongl}
is of the same type as that on the left-hand side of
\eqref{anwatsongl}, this nonterminating $\mathrm A_n$ Watson transformation
does {\em not} reduce to the terminating $\mathrm A_n$ Watson transformation
in \eqref{anwatsongl} as the $\mathrm A_n$ ${}_4\phi_3$ series on the
respective right-hand sides are of different type. Specifically,
\begin{multline}\label{anntwatsongl}
\sum_{k_1,\dots,k_n\ge0}\Bigg(
\prod_{1\le i<j\le n}\frac{x_iq^{k_i}-x_jq^{k_j}}{x_i-x_j}
\prod_{i,j=1}^n\frac{(f_jx_i/x_j;q)_{k_i}}
{(qx_i/x_j;q)_{k_i}}\prod_{i=1}^n\frac{1-ax_iq^{k_i+|k|}}
{1-ax_i}\,\frac{(ax_i;q)_{|k|}}{(ax_iq/f_i;q)_{|k|}}\\\times
\prod_{i=1}^n\frac{(bx_i,cx_i;q)_{k_i}}
{(ax_iq/d,ax_iq/e;q)_{k_i}}\cdot
\frac{(d,e;q)_{|k|}}{(aq/b,aq/c;q)_{|k|}}
\left(\frac{a^2q^2}{bcdef_1\cdots f_n}\right)^{|k|}\Bigg)\\=
\frac{(aq/bf_1\cdots f_n,aq/cf_1\cdots f_n;q)_\infty}{(aq/b,aq/c;q)_\infty}
\prod_{i=1}^n\frac{(ax_iq,af_iq/bcf_1\cdots f_nx_i;q)_\infty}
{(aq/bcf_1\cdots f_nx_i,ax_iq/f_i;q)_\infty}\\\times
\sum_{k_1,\dots,k_n\ge0}\Bigg(
\prod_{1\le i<j\le n}\frac{x_iq^{k_i}-x_jq^{k_j}}{x_i-x_j}
\prod_{i,j=1}^n\frac{(f_jx_i/x_j;q)_{k_i}}
{(qx_i/x_j;q)_{k_i}}\\\times
\prod_{i=1}^n\frac{(ax_iq/de,bx_i,cx_i;q)_{k_i}}
{(ax_iq/d,ax_iq/e,bcf_1\cdots f_nx_i/a;q)_{k_i}}\cdot
q^{|k|}\Bigg)\\
+\frac{(q,a^2q^2/bcdf_1\cdots f_n,a^2q^2/bcef_1\cdots f_n;q)_\infty}
{(a^2q^2/bcdef_1\cdots f_n,aq/b,aq/c;q)_\infty}
\prod_{i=1}^n\frac{(ax_iq;q)_\infty}{(ax_iq/f_i;q)_\infty}\\\times
\sum_{s=1}^n\Bigg[q^{(n-1)k_s}
\frac{(ax_sq/de,bx_s,cx_s;q)_\infty}
{(bcf_1\cdots f_nx_s/aq,ax_sq/d,ax_sq/e;q)_\infty}
\prod_{i=1}^n\frac{(f_ix_s/x_i;q)_\infty}{(qx_s/x_i;q)_\infty}\\\times
\sum_{k_1,\dots,k_n\ge0}\Bigg(
\prod_{\substack{1\le i\le n\\i\neq s}}\frac{x_i}{x_i-x_s}
\prod_{\substack{1\le i<j\le n\\i,j\neq s}}\frac{x_iq^{k_i}-x_jq^{k_j}}{x_i-x_j}
\prod_{\substack{1\le i,j\le n\\i\neq s}}
\frac{(f_jx_i/x_j;q)_{k_i}}{(qx_i/x_j;q)_{k_i}}\\\times
\prod_{\substack{1\le i\le n\\i\neq s}}
\frac{1-bcf_1\cdots f_nx_iq^{k_i-y_s-1}/a}
{1-bcf_1\dots f_nx_i/aq}\,
\frac{(ax_iq/de,bx_i,cx_i;q)_{k_i}}
{(bcf_1\cdots f_nx_i/a,ax_iq/d,ax_iq/e;q)_{k_i}}\\\times
\frac{(a^2q^2/bcdef_1\cdots f_n,aq/bf_1\cdots f_n,aq/cf_1\cdots f_n;q)_{k_s}}
{(q,a^2q^2/bcdf_1\cdots f_n,a^2q^2/bcef_1\cdots f_n;q)_{k_s}}\\\times
\prod_{i=1}^n\frac{(af_iq/bcf_1\cdots f_nx_i;q)_{k_s}}
{(aq^2/bcf_1\cdots f_nx_i;q)_{k_s}}
\cdot\,q^{|k|}\Bigg)\Bigg],
\end{multline}
where $|q|<1$ and $|a^2q^2/bcdef_1\cdots f_n|<1$.

The $f_i=q^{-N_i}$, $i=1,\dots,n$, case of \eqref{anntwatsongl}
gives a terminating $\mathrm A_n$ Watson transformation which is 
different from the one in \eqref{anwatsongl}.
Milne and Newcomb \cite{MN2} obtained yet another nonterminating
$\mathrm A_n$ Watson transformation.

Ito~\cite{I3} derived a $\mathrm{BC}_n$ extension
of Slater's~\cite{Sl} general transformation formula for very-well-poised
balanced ${}_{2r}\psi_{2r}$ series. Some interesting and potentially
useful transformations for $\mathrm A_n$ basic hypergeometric series
involving nested sums were recently given by Fang~\cite{F}.

\smallskip
For further references to summations and transformations
for basic hypergeometric series associated with root systems,
see Milne's survey paper \cite{M9}, 
and Milne and Newcomb's paper \cite{MN2}, and the references therein.

\section{Hypergeometric and basic hypergeometric integrals 
associated with root systems}
\label{sec:int}

There exist a number of hypergeometric integral evaluations
associated with root systems.
Several of them can be viewed as extensions of
Selberg's multivariate extension of 1944 of the classical beta
integral evaluation~\cite{S},
\begin{align}\label{sel}
\int_0^1\cdots\int_0^1
\prod_{1\le i<j\le n}|z_i-z_j|^{2\gamma}\prod_{i=1}^nz_i^{\alpha-1}
(1-z_i)^{\beta-1}\mathrm dz_i&
\notag\\=
\prod_{i=1}^n\frac{\Gamma(\alpha+(i-1)\gamma)\,\Gamma(\beta+(i-1)\gamma)\,
\Gamma(1+i\gamma)}
{\Gamma(\alpha+\beta+(n+i-2)\gamma)\,\Gamma(1+\gamma)}&,
\end{align}
provided $\mathrm{Re}(\alpha)>0$,  $\mathrm{Re}(\beta)>0$ and 
$\mathrm{Re}(\gamma)+\max(\frac 1n,\mathrm{Re}\frac \alpha{n-1},
\mathrm{Re}\frac \beta{n-1})>0$.
The Selberg integral is used in many areas, see 
Chapter~11 of this volume and \cite{FW}.

In 1982, Macdonald~\cite{Mac82} conjectured related constant term
identities associated with root systems together with $q$-analogues.
Assume $R$ to be a reduced root system of rank $n$ with basis of simple
roots $\{\alpha_1,\dots,\alpha_n\}$.
Further, let $e^\alpha$ be the formal exponentials, for $\alpha\in R$,
which form the group ring of the lattice generated by $R$,
and let $d_1,\dots,d_n$ be the degrees of the fundamental invariants
of the Weyl group $W(R)$. Then Macdonald conjectured
\cite[Conjecture~3.1]{Mac82} that for any nonnegative integer $k$
%($k$ can also be $\infty$)
the constant term, i.e.\ the term
not containing any $e^\alpha$, in
\begin{equation}
\prod_{\alpha\in R^+}\prod_{i=1}^k(1-q^{i-1}e^{-\alpha})(1-q^ie^\alpha)
\end{equation}
(where $R^+$ is a system of positive roots in $R$) is
\begin{equation}
\prod_{i=1}^n\frac{(q;q)_{kd_i}}{(q;q)_k\,(q;q)_{k(d_i-1)}}.
\end{equation}
(For the root system $\mathrm A_{n-1}$, this exactly corresponds to the $t=q^k$
case of the squared norm evaluation of Macdonald polynomials indexed by
$\lambda=(0,\dots,0)$, the empty partition, in \eqref{snev}.)
This conjecture can be reformulated in terms of reduced affine
root systems and further strengthened. Generalizations with extra parameter
were proposed by Morris~\cite{M}. A thorough account of the historic
development of $q$-Selberg integrals and corresponding constant term
identities is provided in \cite[Section~2.3]{FW}.

Of particular interest are those multiple integrals which in the univariate
case reduce to the Askey--Wilson integral~\cite{AW}:
\begin{equation}
\frac 1{2\pi \mathrm i}\int_{\mathbb T}
\frac{(z^2,1/z^2;q)_\infty}
{(az,a/z,bz,b/z,cz,c/z,dz,d/z;q)_\infty}
\frac{\mathrm dz}{z}=
\frac{2\,(abcd;q)_\infty}{(q,ab,ac,ad,bc,bd,cd;q)_\infty},
\end{equation}
where $|q|<1$, $|a|<1$, $|b|<1$, $|c|<1$, $|d|<1$,
and $\mathbb T$ is the positively oriented unit circle.
The Askey--Wilson integral is responsible for the orthogonality
of the Askey--Wilson polynomials, which sit at the top
of the $q$-Askey scheme of $q$-orthogonal polynomials.
Such multivariate integral evaluations were first obtained
by Gustafson in the early 1990s. 
In the following, we list some of the Askey--Wilson integral
evaluations associated with root systems.
Many of these or related integrals arise as constant term identities
for (extensions of) Macdonald polynomials. This provides a natural link
of the material presented here with Chapter~9 of this volume.

All these multivariate Askey--Wilson integral evaluations can be further
generalized to multivariate extensions of the Nassrallah--Rahman integral
evaluation \cite[Equation (6.4.1)]{GR}
(which has one more parameter then the Askey--Wilson integral evaluation).
The latter admit elliptic extensions.
They are treated in Section~6.2 of this volume together with some
further extensions to integral transformations.

%Before treating multivariate Askey--Wilson integral evaluations
%which concern integrals on the basic hypergeometric level, we
%review a couple of multivariate hypergeometric integral evaluations.

In the following, let $\mathbb T^n$ be the positively oriented
$n$-dimensional complex torus.
%In the following, let $\mathbb T$ denote 
%$n$-fold direct
%product of the unit circle
%$\{t\in\mathbb C:|t|=1\}$
%traversed in the positive direction.
In \cite[Theorem~6.1]{G7}, the following $\mathrm A_n$
Askey--Wilson integral evaluation was derived:
\begin{align}\label{awint}
\frac 1{(2\pi \mathrm i)^n}\int_{\mathbb T^n}&
\frac{\prod_{1\le i<j\le n+1}(z_i/z_j,z_j/z_i;q)_\infty}
{\prod_{i,j=1}^{n+1}(a_i/z_j,b_iz_j;q)_\infty}
\prod_{i=1}^n\frac{\mathrm dz_i}{z_i}\notag\\
&=\frac{(n+1)!\,\Big(\prod_{i=1}^{n+1}a_ib_i;q\Big)_\infty}{(q;q)_\infty^n\,
\Big(\prod_{i=1}^{n+1}a_i,\prod_{i=1}^{n+1}b_i;q\Big)_\infty
\prod_{i,j=1}^{n+1}(a_ib_j;q)_\infty},
\end{align}
where $\prod_{i=1}^{n+1}z_i=1$, provided $|q|<1$, $|a_i|<1$ and $|b_i|<1$,
for $1\le i\le n+1$.

A considerably more complicated $\mathrm A_n$
Askey--Wilson integral evaluation, depending
on the parity of $n$, was given in \cite[Theorem~1.1]{GuR}:
\begin{multline}
\frac 1{(2\pi \mathrm i)^n}\int_{\mathbb T^n}
\frac{\prod_{1\le i<j\le n+1}(z_i/z_j,z_j/z_i;q)_\infty}
{\prod_{i,j=1}^{n+1}(a_i/z_j;q)_\infty
\prod_{1\le i<j\le n+1}(bz_iz_j;q)_\infty}
\prod_{i=1}^{n+1}\frac{(S/z_i;q)_\infty}
{\prod_{j=1}^3(bc_jz_i;q)_\infty}
\frac{\mathrm dz_i}{z_i}\\
=\begin{cases}\displaystyle{
\frac{(n+1)!\,\Big(b^{(n+4)/2}\prod_{i=1}^{n+1}a_i
\prod_{j=1}^3c_j;q\Big)_\infty\prod_{i=1}^{n+1}(S/a_i;q)_\infty}
{(q;q)_\infty^n\,\Big(\prod_{i=1}^{n+1}a_i,b^{(n+4)/2}\prod_{j=1}^3c_j;q\Big)_\infty
\prod_{j=1}^3(b^{(m+2)/2}c_j;q)_\infty}}&{}\\[1em]
\displaystyle{\times\;
\frac{\prod_{j=1}^3\Big(b^{(n+2)/2}c_j\prod_{i=1}^{n+1}a_i;q\Big)_\infty}
{\prod_{i=1}^{n+1}\prod_{j=1}^3(ba_ic_j;q)_\infty
\prod_{1\le i<j\le n+1}(ba_ia_j;q)_\infty},}
&\text{for $n$ even},\\[2em]
\displaystyle{\frac{(n+1)!\,\Big(b^{(n+1)/2}\prod_{i=1}^{n+1}a_i;q\Big)_\infty
\prod_{i=1}^{n+1}(S/a_i;q)_\infty}
{(q;q)_\infty^n\,\big(b^{(n+1)/2},\prod_{i=1}^{n+1}a_i;q\big)_\infty
\prod_{i=1}^{n+1}\prod_{j=1}^3(ba_ic_j;q)_\infty}}&{}\\[1em]
\displaystyle{\times\;
\frac{\prod_{j=1}^3\Big(b^{(n+3)/2}
\prod_{i=1}^{n+1}a_i\prod_{\substack{1\le k\le 3\\k\ne j}}c_k;q\Big)_\infty}
{\prod_{1\le i<j\le 3}(b^{(n+3)/2}c_ic_j;q)_\infty
\prod_{1\le i<j\le n+1}(ba_ia_j;q)_\infty},}
&\text{for $n$ odd}.
\end{cases}
\end{multline}
where $\prod_{i=1}^{n+1}z_i=1$ and
$S=b^{n+2}\prod_{i=1}^{n+1}a_i\prod_{j=1}^3c_j$, provided $|q|<1$, $|b|<1$,
$|a_i|<1$ and $|c_j|<1$, for $1\le i\le n+1$ and $1\le j\le n$ .

The following $\mathrm C_n$
Askey--Wilson integral evaluation was derived in \cite[Theorem 7.1]{G7}:
\begin{align}
\frac 1{(2\pi \mathrm i)^n}\int_{\mathbb T^n}
\frac{\prod_{1\le i<j\le n}
(z_i/z_j,z_j/z_i,z_iz_j,1/z_iz_j;q)_\infty}
{\prod_{i=1}^{2n+2}\prod_{j=1}^n(a_iz_j,a_i/z_j;q)_\infty}
\prod_{i=1}^n(z_i^2,1/z_i^2;q)_\infty\frac{\mathrm dz_i}{z_i}\notag&\\
=\frac{2^nn!\,\Big(\prod_{i=1}^{2n+2}a_i;q\Big)_\infty}{(q;q)_\infty^n
\prod_{1\le i<j\le 2n+2}(a_ia_j;q)_\infty}&,
\end{align}
provided $|q|<1$ and $|a_i|<1$ for $1\le i\le n$.
%provided  $a_ia_j\neq q^k$, for $1\le i\le j\le 2n+2$ and $k=0,-1,-2,\dots$,
%where the contour $C$ is the unit circle traversed in the positive direction
%but with suitable deformations to separate the sequences of poles converging
%to zero from the sequences of poles diverging to infinity.

Another $\mathrm C_n$ Askey--Wilson integral evaluation was given in
\cite[Equation (2)]{G5}:
\begin{align}\label{intmk}
\frac 1{(2\pi \mathrm i)^n}\int_{\mathbb T^n}\prod_{1\le i<j\le n}
\frac{(z_i/z_j,z_j/z_i,z_iz_j,1/z_iz_j;q)_\infty}
{(bz_i/z_j,bz_j/z_i,bz_iz_j,b/z_iz_j;q)_\infty}
\prod_{i=1}^n\frac{(z_i^2,1/z_i^2;q)_\infty}
{\prod_{j=1}^4(a_jz_i,a_j/z_j;q)_\infty}\frac{\mathrm dz_i}{z_i}\notag&\\
=\frac{2^nn!\,(b;q)_\infty^n}{(q;q)_\infty^n}
\prod_{i=1}^n\frac{\Big(b^{n+i-2}\prod_{j=1}^4a_j;q\Big)_\infty}
{(b^i;q)_\infty\prod_{1\le j<k\le 4}(a_ja_kb^{i-1};q)_\infty}&,
\end{align}
provided $|q|<1$, $|b_i|<1$ and $|a_j|<1$, for $1\le i\le n$ and $1\le j\le 4$.
By suitably specializing the variables $a_j$ for $1\le j\le 4$,
this multivariate integral evaluation can be used (see \cite{G5}) to
prove Morris' \cite{M} generalizations of the
Macdonald conjectures \cite{M}
for the affine root systems
$\mathrm C_n$, $\mathrm C_n^\vee$, $\mathrm{BC}_n$, $\mathrm B_n$
$\mathrm B_n^\vee$ and $\mathrm D_n$ (using Macdonald's
classification in \cite{Mac72}).

The multivariate integral evaluation in \eqref{intmk} explicitly describes
the normalization factor of the orthogonality measure for the
Mac\-donald--Koornwinder polynomials (see \cite{Ko} and Chapter 9
of this volume), the $\mathrm{BC}_n$ generalization of the
Askey--Wilson polynomials. 

In \cite[Theorem~8.1]{G7} an Askey--Wilson integral evaluation
for the root system $\mathrm G_2$ was given:
\begin{align}
\frac 1{(2\pi \mathrm i)^2}\int_{\mathbb T^2}&
\frac{\prod_{\substack{1\le i,j\le 3\\i\ne j}}(z_i/z_j;q)_\infty
\prod_{j=1}^3(z_j,1/z_j;q)_\infty}
{\prod_{i=1}^4\prod_{j=1}^3(a_iz_j,a_i/z_j;q)_\infty}
\frac{\mathrm dz_1}{z_1}\frac{\mathrm dz_2}{z_2}\notag\\
&=\frac{12\;\Big(\prod_{i=1}^4a_i^2;q\Big)_\infty\prod_{i=1}^4(a_i;q)_\infty}
{(q;q)_\infty^2\,\Big(\prod_{i=1}^4a_i;q\Big)_\infty
\prod_{1\le i\le j\le 4}(a_ia_j;q)_\infty
\prod_{1\le i<j<k\le 4}(a_ia_ja_k;q)_\infty},
\end{align}
where $\prod_{j=1}^3z_j=1$ and $|a_i|<1$ for $1\le i\le 4$.

All these basic hypergeometric integral evaluations can be specialized
to ordinary hypergeometric integral evaluations by taking suitable limits.
In particular, if in \eqref{awint} one replaces the parameters
$a_i$ by $q^{a_i}$, and $b_i$ by $q^{b_i}$, for $1\le i\le n+1$,
and then takes the limit as $q\to 1^-$, one obtains the following
multidimensional Mellin--Barnes integral
\cite[Theorem~9.1]{G5}
\begin{multline}\label{awinto}
\frac 1{(2\pi\mathrm i)^n}\int_{-\mathrm{i}\infty}^{\mathrm{i}\infty}
\cdots\int_{-\mathrm{i}\infty}^{\mathrm{i}\infty}
\frac{\prod_{i,j=1}^{n+1}\Gamma(a_i-z_j)\Gamma(b_i+z_j)}
{\prod_{\substack{1\le i,j\le n+1\\i\ne j}}\Gamma(z_i-z_j)}
\prod_{i=1}^n\frac{\mathrm dz_i}{z_i}\\
=\frac{(n+1)!\,
\Gamma(a_1+\cdots+a_{n+1})\,\Gamma(b_1+\cdots+b_{n+1})
\prod_{i,j=1}^{n+1}\Gamma(a_i+b_j)}{
\Gamma(a_1+\cdots+a_{n+1}+b_1+\cdots+b_{n+1})},
\end{multline}
where $\sum_{i=1}^{n+1}z_i=0$, provided $\mathrm{Re}(a_i)>0$ and
$\mathrm{Re}(b_i)>0$,
for $1\le i\le n+1$. (For a generalization of \eqref{awinto},
obtained by taking a suitable $q\to 1^-$ limit from an $\mathrm A_n$
Nassrallah--Rahman integral that extends \eqref{awint}, see
\cite[Theorem~5.1]{G6}).

Here we reproduced only a few of the many existing integral evaluations.
More can be found in the papers \cite{G7,G8,Ra1,SW}
(to list just a few relevant references).
An interesting integral transformation with $\mathrm F_4$ symmetry
has been given in \cite{vdB}.
For further discussion of integral identities
(evaluations and transformations) associated with root systems,
where such identities are considered at the elliptic level,
see Section~6.2 of this volume.

\section{Basic hypergeometric series with Macdonald
polynomial argument}\label{sec:macd}

The series considered here were first introduced by Macdonald
(in unpublished work \cite{Mac1} of 1987), and by Kaneko~\cite{Ka1}.

Important special cases were considered earlier.
Basic hypergeometric series with Schur polynomial argument
(the Schur polynomials correspond to the $q=t$ case of the
Macdonald polynomials) were studied by Milne~\cite{M92}
who derived  ${}_1\phi_0$, ${}_2\phi_1$ and ${}_1\psi_1$ summations
and several transformations for such series.
Hypergeometric series with Jack polynomial argument (the Jack
polynomials indexed by $\alpha$ correspond to the $q=t^{\alpha}$, $t\to 1$
specialization of the Macdonald polynomials) were studied by
Herz, Constantine and Muirhead~\cite{Con,He,Mu} for $\alpha=2$
(the zonal polynomial case) and for arbitrary $\alpha$ by
Kor\'anyi, Yan and Kaneko~\cite{Ka0,Kor,Y1,Y2}.

For a thorough treatment of Macdonald polynomials (by which we mean
the $\mathrm{GL}_n$ type symmetric Macdonald polynomials in the terminology
of Chapter 9 of this volume; for the general root system case see
\cite{Mac2} and Chapter 9 of this volume),
see \cite[Chapter~VI]{Mac0} and Sections~9.1.1 and 9.3.7 of this volume.
(Macdonald's book \cite{Mac0} also
deals thoroughly with important special cases of the
Macdonald polynomials, including in particular the aforementioned
Schur, zonal and Jack polynomials.)
See Chapter~10 of this volume for a survey on combinatorial aspects
of these multivariate polynomials.

Let $\Lambda_{n}$ denote the ring of symmetric functions
in the variables $z=(z_1,\dots,z_n)$ over $\mathbb C$.
Further, we assume two nonzero generic parameters $q,t$ satisfying
$|q|,|t|<1$.
The Macdonald polynomials $P_\lambda(z;q,t)$ (often shortened to $P_\lambda$
or $P_\lambda(z)$ as long as no ambiguity arises), indexed by partitions
$\lambda$ of length $l(\lambda)\le n$,
form an orthogonal basis of $\Lambda_{n}$.
They can be defined as the unique family of symmetric polynomials
whose expansion in terms of the monomial symmetric functions
$m_\lambda(z)$ is uni-upper-triangular with respect to the dominance
order $<$ of partitions,
\begin{equation*}
P_{\lambda}(z;q,t)=m_\lambda(z)+\sum_{\mu<\lambda}c_{\lambda\mu}(q,t)m_\mu(z)
\end{equation*}
(with $c_{\lambda\mu}(q,t)$ being a rational function in $q$ and $t$),
being orthogonal with respect to the scalar product
\cite[Chapter~VI, Section~9]{Mac0}
\begin{equation*}
\langle f,g\rangle=\frac 1{n!(2\pi\mathrm i)^n}\int_{\mathbb T^n}
f(z)\overline{g(z)}\,\Delta_{q,t}(z)
\prod_{i=1}^n\frac{\mathrm dz_i}{z_i},\qquad \text{for $f,g\in\Lambda_{n}$}, 
\end{equation*}
where
\begin{equation*}
\Delta_{q,t}(z)=\prod_{\substack{1\le i,j\le n\\i\neq j}}
\frac{(z_i/z_j;q)_\infty}{(tz_i/z_j;q)_\infty}.
\end{equation*}
As in Section~\ref{sec:int}, $\mathbb T^n$ is the positively oriented
$n$-dimensional complex torus.
The squared norm evaluation of $P_\lambda$ is
\cite[Chapter~VI, Section~9, Example~1.(d)]{Mac0}
\begin{equation}\label{snev}
\langle
P_{\lambda},P_{\lambda}
\rangle
=\prod_{1\le i<j\le n}\frac{(q^{\lambda_i-\lambda_j}t^{j-i},
q^{\lambda_i-\lambda_j+1}t^{j-i};q)_\infty}{(q^{\lambda_i-\lambda_j}t^{j-i+1},
q^{\lambda_i-\lambda_j+1}t^{j-i-1};q)_\infty}.
\end{equation}

In \cite[Chapter~VI]{Mac0} Macdonald develops most of
the theory for the polynomials $P_{\lambda}(z;q,t)$ using
a different (albeit, up to normalization, equivalent)
scalar product (which we are not displaying here)
that is more algebraic in nature and does not require the
conditions $|q|<1$ and $|t|<1$. Rather than considering
symmetric functions over $\mathbb C$,
Macdonald assumes $q$ and $t$ to be indeterminate and considers
symmetric functions over $\mathbb Q(q,t)$.
The above scalar product has the advantage that the structure
of the root system $\mathrm A_{n-1}$ is clearly visible. This aspect
of the theory generalizes to other root systems, see \cite{Mac2}
and Chapter 9 of this volume.

The $P_\lambda(z;q,t)$ are homogeneous in $z=(z_1,\dots,z_n)$
of degree $|\lambda|$.
They satisfy the stability property
\begin{equation*}
P_{\lambda}(z_1,\dots,z_n;q,t)=P_{\lambda}(z_1,\dots,z_n,0;q,t).
\end{equation*}
%for $l(\lambda)\le n$,
%viewed as an identity in the larger ring
%$\Lambda_{n+1,F}\supset\Lambda_{n,F}$.
Further, they satisfy \cite[Chapter~VI, Equation~(4.17)]{Mac0}
\begin{equation}\label{shift}
P_{\lambda}(z;q,t)=(z_1\cdots z_n)^{\lambda_n}P_{\lambda-\lambda_n}(z;q,t),
\end{equation}
where
$\lambda-\lambda_n:=(\lambda_1-\lambda_n,\dots,\lambda_{n-1}-\lambda_n,0)$
for any partition $\lambda$ with $l(\lambda)\le n$.

For any partition $\lambda$, and $f\in\Lambda_{n}$
let $u_\lambda:\Lambda_{n}\to \mathbb C$ be the
evaluation homomorphism defined by 
\begin{equation*}
u_\lambda\big(f(z)\big)=\begin{cases}
f(z)\Bigr\rvert_{z_i=q^{\lambda_i} t^{n-i},\,1\le i\le n}&
\text{for $l(\lambda)\le n$},\\
0& \text{otherwise}.
\end{cases}
\end{equation*}
The following evaluation symmetry (cf.\
\cite[Chapter~VI, Equation (6.6)]{Mac0}), first proved by Koornwinder in
unpublished work, is very useful (in particular, for interchanging summations
in the process of deriving transformations):
\begin{equation}\label{evsym}
u_0(P_{\lambda})\,u_{\lambda}(P_{\mu})
=u_0(P_{\mu})\,u_{\mu}(P_{\lambda}).
\end{equation}

For any partition $\lambda$, let
\begin{equation*}
(a;q,t)_\lambda=\prod_{i\ge 1}(at^{1-i};q)_{\lambda_i},
\end{equation*}
and
\begin{equation*}
n(\lambda)=\sum_{i\ge 1}(i-1)\lambda_i=\sum_{i\ge 1}\binom{\lambda_i'}2,
\end{equation*}
where $\lambda'$ denotes the conjugate partition of $\lambda$.
We also use the shorthand notation
\begin{equation*}
(a_1,\dots,a_n;q,t)_\lambda=(a_1;q,t)_\lambda\cdots(a_n;q,t)_\lambda.
\end{equation*}
Further, for $l(\lambda)\le n$, we define 
\begin{subequations}
\begin{align}
c_\lambda(q,t)&=\prod_{i=1}^n(t^{n-i+1};q)_{\lambda_i}\prod_{1\le i<j\le n}
\frac{(t^{j-i};q)_{\lambda_i-\lambda_j}}{(t^{j-i+1};q)_{\lambda_i-\lambda_j}},\\
c'_\lambda(q,t)&=\prod_{i=1}^n(qt^{n-i};q)_{\lambda_i}\prod_{1\le i<j\le n}
\frac{(qt^{j-i-1};q)_{\lambda_i-\lambda_j}}{(qt^{j-i};q)_{\lambda_i-\lambda_j}}.
\end{align}
\end{subequations}

An important normalization of the Macdonald polynomials is given by
\begin{equation*}
Q_{\lambda}(z;q,t)=b_{\lambda}(q,t)P_{\lambda}(z;q,t),
\end{equation*}
where
\begin{equation*}
b_{\lambda}(q,t)=\frac{h_\lambda(q,t)}{c'_\lambda(q,t)}.
\end{equation*}
The $Q_\lambda$ are exactly the polynomials dual to $P_\lambda$
with respect to scalar product mentioned right after Equation~\eqref{snev}.
The two normalizations of Macdonald polynomials appear jointly
in the Cauchy identity
\begin{equation}\label{cauchy}
\sum_{\lambda}P_{\lambda}(z;q,t)Q_{\lambda}(y;q,t)=\prod_{i,j=1}^n
\frac{(tz_iy_j;q)_\infty}{(z_iy_j;q)_\infty}.
\end{equation}

Let $a$ be an indeterminate and define the homomorphism
\begin{equation*}
\epsilon_{a;t}:\Lambda_{n}\to \mathbb C[a]
\end{equation*}
by its action on the power sum symmetric functions $p_r=p_r(z_1,\dots,z_n):=
\sum_{i=1}^nz_i^r$ for $r\ge 1$
(which algebraically generate $\Lambda_{n}$), namely
\begin{equation*}
\epsilon_{a;t}(p_r)=\frac{1-a^r}{1-t^r},
\end{equation*}
for each $r\ge 1$. For $a=t^n$, we have
$\epsilon_{t^n;t}(f)=f(1,t,\dots,t^{n-1})$ for any $f\in\Lambda_{n}$.

The following evaluations are useful (cf.\ \cite[Theorem~3.3]{Ka1}):
\begin{equation}\label{specmac}
\epsilon_{a;t}\big(P_\lambda(z;q,t)\big)=t^{n(\lambda)}
\frac{(a;q,t)_{\lambda}}{c_\lambda(q,t)},
\qquad
\epsilon_{a;t}\big(Q_\lambda(z;q,t)\big)=t^{n(\lambda)}
\frac{(a;q,t)_{\lambda}}{c'_\lambda(q,t)}.
\end{equation}

Basic hypergeometric series with Macdonald polynomial argument are defined as
\begin{equation}
_r\Phi_s\left[\begin{matrix}a_1,\dots,a_r\\
b_1,\dots,b_s\end{matrix};q,t,z\right]
=\sum_\lambda\left((-1)^{|\lambda|}q^{n(\lambda')}t^{-n(\lambda)}\right)^{s+1-r}
\frac{t^{n(\lambda)}}{c'_{\lambda}(q,t)}
\frac{(a_1,\dots,a_r;q,t)_\lambda}
{(b_1,\dots,b_s;q,t)_\lambda}P_\lambda(z;q,t),
\end{equation}
provided that the series converges.
%\begin{align}
%&_r\Phi_s\left[\begin{matrix}a_1,\dots,a_r\\
%b_1,\dots,b_s\end{matrix};q,t,z\right]\notag\\
%&=\sum_\lambda\left((-1)^{|\lambda|}
%q^{n(\lambda')}t^{-n(\lambda)}\right)^{s+1-r}
%\frac{t^{n(\lambda)}}{c'_{\lambda}(q,t)}
%\frac{(a_1,\dots,a_r;q,t)_\lambda}
%{(b_1,\dots,b_s;q,t)_\lambda}P_\lambda(z;q,t).
%\end{align}

Application of the homomorphism $\epsilon_{a;t}$ with respect to $y$
to both sides of the Cauchy identity in \eqref{cauchy}
immediately gives the following $q$-binomial theorem for Macdonald polynomials:
\begin{equation}\label{qbinmac}
_1\Phi_0\left[\begin{matrix}a\\-\end{matrix};q,t,z\right]=
\prod_{i=1}^n\frac{(az_i;q)_\infty}{(z_i;q)_\infty},
\end{equation}
which converges for $|z_i|<1$, $1\le i\le n$. It is interesting that
the right-hand side is independent of $t$.

Baker and Forrester~\cite{BF} made use of the
$q$-binomial theorem for Macdonald polynomials and the
evaluation symmetry \eqref{evsym} to derive
the following multivariate generalization of the Heine transformation:
\begin{equation}\label{heine}
_2\Phi_1\left[\begin{matrix}a,b\\c\end{matrix};q,t,xt^\delta\right]=
\prod_{i=1}^n\frac{(bt^{1-i},axt^{n-i};q)_\infty}{(ct^{1-i},xt^{n-i};q)_\infty}\cdot
{}_2\Phi_1\left[\begin{matrix}c/b,xt^{n-1}\\axt^{n-1}\end{matrix};
q,t,bt^{1-n}t^\delta\right],
\end{equation}
valid for $|x|<1$ and $|bt^{1-n}|<1$,
where $xt^\delta$ stands for the argument $(x,xt,\dots,xt^{n-1})$.
Notice that this transformation involves specialized Macdonald polynomials
(which factorize since
$P_\lambda(xt^\delta)=x^{|\lambda|}t^{n(\lambda)}(t^n;q)_\lambda/c_\lambda(q,t)$
due to homogeneity, and the specialization \eqref{specmac}) on both sides.
A multivariate generalization of the first iterate of the Heine transformation
involving unspecialized interpolation Macdonald polynomials was given by
Lascoux, Rains and Warnaar~\cite[Corollary~10.2]{LRW}. A further
extension was obtained by Lascoux and Warnaar~\cite[Corollary~6.3]{LW}
as a special case of of multivariate extension of the
$q$-Kummer--Thomae--Whipple transformation \cite[Corollary~6.2]{LW}.

For $z=ct^{1-n}/ab$ the right-hand side of \eqref{heine} reduces to
a\/ ${}_1\Phi_0$ series which can be summed using \eqref{qbinmac}.
This gives a multivariate extension of the $q$-Gau{\ss} summation:
\begin{equation}\label{qgauss}
_2\Phi_1\left[\begin{matrix}a,b\\c\end{matrix};q,t,
\frac {ct^{1-n}}{ab}t^\delta\right]=
\prod_{i=1}^n\frac{(ct^{1-i}/b,ct^{1-i}/a;q)_\infty}{(ct^{1-i}/ab,ct^{1-i};q)_\infty},
\end{equation}
valid for $|ct^{1-n}/ab|$.
More general $q$-Gau{\ss} summations involving unspecialized
(non-)sym\-met\-ric Macdonald polynomials were given by Lascoux, Rains and
Warnaar~\cite{LRW}, and by Lascoux and Warnaar \cite[Corollary~5.4]{LW}.

For general unspecialized argument $z=(z_1,\dots,z_n)$
Baker and Forrester~\cite{BF}, building on work of Kaneko~\cite{Ka1},
proved the following multivariate extension of the Euler transformation
(or equivalently, of the second iterate of the Heine transformation):
\begin{equation}\label{qeuler}
_2\Phi_1\left[\begin{matrix}a,b\\c\end{matrix};q,t,z\right]=
\prod_{i=1}^n\frac{(abz_i/c;q)_\infty}{(z_i;q)_\infty}\cdot
{}_2\Phi_1\left[\begin{matrix}c/a,c/b\\c\end{matrix};
q,t,abz/c\right],
\end{equation}
valid for $|z_i|<1$ and $|abz_i/c|<1$, $1\le i\le n$.
A nonsymmetric extension was given in \cite[Corollary~10.3]{LRW}.

We list two other results from \cite{BF}. Let $N$ be a nonnegative integer.
The $q$-Pfaff--Saalsch\"utz summation for basic hypergeometric
series with specialized Macdonald polynomial argument is
\begin{equation}\label{qps}
_3\Phi_2\left[\begin{matrix}a,b,q^{-N}\\c,abq^{1-N}t^{n-1}/c\end{matrix};
q,t,qt^\delta\right]=\prod_{i=1}^n
\frac{(ct^{1-i}/a,ct^{1-i}/b;q)_N}{(ct^{1-i},ct^{1-i}/ab;q)_N}.
\end{equation}
This can be generalized to a multivariate Sears' transformation
with specialized Macdonald polynomial arguments:
\begin{align}\label{sears}
&{}_4\Phi_3\left[\begin{matrix}a,b,c,q^{-N}\\d,e,ft^{n-1}\end{matrix};
q,t,qt^\delta\right]\notag\\
&=a^{nN}
\prod_{i=1}^n\frac{(et^{1-i}/a,ft^{n-i}/a;q)_N}{(et^{1-i},ft^{n-i};q)_N}\cdot
{}_4\Phi_3\left[\begin{matrix}a,d/b,d/c,q^{-N}\\
d,aq^{1-N}t^{n-1}/e,aq^{1-N}/f\end{matrix};
q,t,qt^\delta\right],
\end{align}
where $def=abcq^{1-N}$.
In \cite[Section~4]{Ra0} Rains proved extensions of \eqref{qps} and
\eqref{sears} for Macdonald polynomials indexed by partitions of skew shape.
Extensions of \eqref{qps} and \eqref{sears} 
to nonsymmetric Macdonald polynomials
were given in \cite[Theorem~6.6 and Proposition~6.8]{LRW}.

Kaneko~\cite{Ka1} developed $q$-difference equations
for basic hypergeometric series with Macdonald polynomial argument
and related them to $q$-Selberg integrals.
Warnaar~\cite{W3} proved various generalizations of
$q$-Selberg integral evaluations and constant term identities,
including a $q$-analogue of the Hua--Kadell formula for Jack polynomials
(cf.\ \cite[Theorem~5.2.1]{Hu} and \cite[Theorem~2]{Kad}).
Rains and Warnaar~\cite[Section~5.3]{RW} obtained further multivariate
${}_4\Phi_3$ transformations.

We complete this section with a multivariate extension of
Ramanujan's $_1\psi_1$ summation formula due to Kaneko~\cite{Ka2},
which is a $t$-extension of an earlier result by Milne~\cite{M92}
(namely, for basic hypergeometric series with Schur function argument).

Let $\mathbb Z^n_{\ge}=\{(\lambda_1,\dots,\lambda_n)\in\mathbb Z^n\mid
\lambda_1\ge\dots\ge\lambda_n\}$. By \eqref{shift}, the Macdonald
polynomials $P_\lambda(z;q,t)$ can be defined for any
$\lambda\in\mathbb Z^n_{\ge}$.
Bilateral basic hypergeometric series with Macdonald polynomial argument
are defined as
\begin{align}
_r\Psi_{s+1}\left[\begin{matrix}a_1,\dots,a_r\\
b,b_1,\dots,b_s\end{matrix};q,t,z\right]\notag
=\sum_{\lambda\in\mathbb Z^n_{\ge}}\Bigg(&
\left((-1)^{|\lambda|}q^{n(\lambda')}t^{-n(\lambda)}\right)^{s+1-r}
\frac{(q;q)_\infty^n}{(b;q)_\infty^n}\prod_{i=1}^n
\frac{(bq^{\lambda_i}t^{n-i};q)_\infty}{(q^{\lambda_i+1}t^{n-i};q)_\infty}\\
%\notag\\&\qquad\qquad\qquad\qquad
&\times
\frac{t^{n(\lambda)}}{c'_{\lambda}(q,t)}
\frac{(a_1,\dots,a_r;q,t)_\lambda}
{(b_1,\dots,b_s;q,t)_\lambda}
P_\lambda(z;q,t)\Bigg),
\end{align}
provided that the series converges.
With this notation, Kaneko's $_1\psi_1$ summation for Macdonald polynomials is
\begin{equation}\label{kan1psi1}
_1\Psi_{1}\left[\begin{matrix}a\\b\end{matrix};q,t,z\right]=\prod_{i=1}^n
\frac{(qt^{n-i},bt^{1-i}/a,az_i,q/az_i;q)_\infty}
{(bqt^{1-i},qt^{n-i}/a,z_i,bt^{1-n}/az_i;q)_\infty},
\end{equation}
subject to $|bt^{1-n}/a|<|z_i|<1$, for $1\le i\le n$.

In \cite[Theorem~2.6]{W5}, Warnaar gives a generalization of
\eqref{kan1psi1} involving a pair of Macdonald polynomials
in two independent sets of variables. Other identities of this
type are obtained in \cite{W4}.

\section{Remarks on applications}\label{sec:appl}

As mentioned in the introduction, hypergeometric series
associated to root systems first arose in the context of
$3j$ and $6j$ symbols for the unitary groups
\cite{AJJ,CCB,HBL}. This initiated their study and that of their
basic analogues from a pure mathematics point of view.

Basic hypergeometric series associated with root systems
have found applications in various areas. We list a few
occurrences, making no claim about completeness.
First of all, such series, in particular
multivariate ${}_6\psi_6$ summations associated with root systems,
were used to give elementary proofs of the
Macdonald identities~\cite{G3,M2}.
More generally, these series were used for
deriving expansions of various special powers of the eta
function~\cite{BW,LeMa,LeMb,M00,WZ} and for
establishing infinite families of exact formulae for sums of
squares and of triangular numbers \cite{M8,R6,R7}.
Basic hypergeometric series associated with root systems
were also employed in the enumeration of plane partitions~\cite{GK,KS,R8}.
Applications to Macdonald polynomials were given in \cite{KN,Schl4}.
Basic hypergeometric integrals associated to root systems were used
in the construction of $\mathrm{BC}_n$ orthogonal
polynomials and $\mathrm{BC}_n$ biorthogonal
rational functions that generalize the Macdonald polynomials,
see \cite{Ko,Ra0}.
Watson transformations (and related transformations)
associated to root systems were used in~\cite{Ba,BW,C,GOW,ZW} to derive
multiple Rogers--Ramanujan identities and characters for affine Lie algebras.
For applications to quantum groups, see~\cite{R2}.
Basic hypergeometric series of Macdonald polynomial argument
were used to construct Selberg-type integrals for $\mathrm A_{n-1}$~\cite{W6}.
Also, hypergeometric series with Jack and zonal polynomial argument
appeared in studies on random matrices~\cite{DL,FR} and
Selberg integrals~\cite{Ka0,Kor}.

Very recently the subject has gained growing attention by physicists
working in spin models and in quantum field theory.
In particular, it was shown in \cite{DM,DMV} that
Gustafson's multivariate hypergeometric integrals appear naturally in the
integrable spin models.
Further, it was shown \cite{HHL,J}
that the partition functions in $3d$ field theories
can be expressed in terms of specific basic hypergeometric integrals.
As made explicit in \cite{DSV}, these partition functions can also be obtained
by reduction from the $4d$ superconformal indices which,
according to \cite{DO},
can be identified with elliptic hypergeometric integrals.
(The latter are reviewed in Chapter~6 of this volume.)
Accordingly, multivariate basic hypergeometric integrals and series
associated with various symmetry groups (or gauge groups, in the terminology
of quantum field theory) appear as explicit expressions for the respective
partition functions~\cite{DSV,SV1,SV2}. Several of these are new and
await further mathematical study.
These partition functions can also be interpreted as solutions of the
Yang--Baxter equation, see \cite{G,GS}.

\end{document}